\documentclass[11pt,reqno,a4paper]{amsart}

\usepackage[T1]{fontenc}
\usepackage[utf8]{inputenc}
\usepackage{lmodern}
\usepackage[english]{babel}

\usepackage[
colorlinks=true,
urlcolor=teal,
citecolor=magenta,
linkcolor=teal,
breaklinks=true]
{hyperref}
\usepackage{amsthm}
\usepackage{thmtools}
\usepackage[noabbrev,capitalize,nameinlink]{cleveref}

\usepackage{amsmath}
\usepackage{amssymb}

\usepackage{mathtools}
\usepackage{mathrsfs}
\usepackage[normalem]{ulem}

\usepackage{bookmark}

\usepackage[initials,msc-links,
]{amsrefs}

\usepackage[dvipsnames]{xcolor}
\usepackage{graphicx}
\usepackage{url}

\usepackage{enumitem}

\usepackage{aurical} 
\usepackage{changes}

\usepackage{tikz-cd}
\usepackage{float}
\usepackage{enumitem}

\theoremstyle{plain}
\newtheorem{theorem}{Theorem}[section]
\newtheorem{lemma}[theorem]{Lemma}
\newtheorem{proposition}[theorem]{Proposition}
\newtheorem{corollary}[theorem]{Corollary}

\theoremstyle{definition}
\newtheorem{definition}[theorem]{Definition}

\newtheorem{remark}[theorem]{Remark}

\definecolor{deepgreen}{rgb}{0.1,0.6,0.1}

\newcommand{\g}{\mathfrak g}

\newcommand{\N}{\mathbb{N}}

\newcommand{\R}{\mathbb{R}}

\newcommand{\average}{{\mathchoice {\kern1ex\vcenter{\hrule height.4pt
    width 6pt
    depth0pt} \kern-9.7pt} {\kern1ex\vcenter{\hrule height.4pt width 4.3pt
    depth0pt}
   \kern-7pt} {} {} }}

\newcommand{\E}{\mathrm{E}}
\newcommand{\End}{\mathrm{End}}

\DeclareMathOperator{\ad}{ad}

\DeclareMathOperator{\Ad}{Ad}

\allowdisplaybreaks

\newcommand{\df}{\mathrm{d}}

\newcommand{\nomedascegliere}{{star}}

\newcommand{\starr}{\mathrm{star}}
\newcommand{\interior}{\operatorname{int}}

\definechangesauthor[name=Nicola, color=brown]{N}
\definechangesauthor[name=Sebastiano, color=blue]{S}
\definechangesauthor[name=Enrico, color=magenta]{E}

\begin{document}

\title[Normal curves in sub-Finsler Lie groups]{Normal curves in sub-Finsler Lie groups: Branching for strongly convex norms and face stability for polyhedral norms}

 \author[E. Le Donne]{Enrico Le Donne}
\address[E. Le Donne]{University of Fribourg, Chemin du Mus\'ee~23, 1700 Fribourg, Switzerland}
\email{\href{mailto:enrico.ledonne@unifr.ch}{enrico.ledonne@unifr.ch}}
 \author[S. Nicolussi Golo]{Sebastiano Nicolussi Golo}
 \address[S. Nicolussi Golo]{University of Fribourg, Chemin du Mus\'ee~23, 1700 Fribourg, Switzerland}
\email{\href{mailto:sebastiano.nicolussigolo@unifr.ch}{sebastiano.nicolussigolo@unifr.ch}}
 \author[N. Paddeu]{Nicola Paddeu}
 \address[N. Paddeu]{University of Fribourg, Chemin du Mus\'ee~23, 1700 Fribourg, Switzerland  
}
\email{\href{mailto:nicola.paddeu@unifr.ch}{nicola.paddeu@unifr.ch}}

\date{\today}

\keywords{Carnot groups, sub-Finsler, geodesics, normal curves}


\thanks{\textit{Memberships and funding information.} E. Le Donne, S. Nicolussi Golo, and N. Paddeu were partially supported by the Swiss National Science Foundation
 (grant 200021-204501 `\emph{Regularity of sub-Riemannian geodesics and applications}'). }

\begin{abstract}
We consider Lie groups equipped with left-invariant subbundles of their tangent bundles and norms on them. On these sub-Finsler structures, we study the normal curves in the sense of control theory. We revisit the Pontryagin Maximum Principle using tools from convex analysis, expressing the normal equation as a differential inclusion involving the subdifferential of the dual norm. In addition to several properties of normal curves, we discuss their existence, the possibility of branching, and local optimality. Finally, we focus on polyhedral norms and show that normal curves have controls that locally take values in a single face of a sphere with respect to the norm.
\end{abstract}
\vspace*{-1.0cm}
\bibliographystyle{abbrv}
\maketitle

\tableofcontents

\vspace{.5cm}

%
%
%



\section{Introduction}

The study of geodesics in sub-Finsler Lie groups lies at the intersection of two important areas of research in mathematics: Control Theory and Geometric Group Theory. On the one hand, geodesics in sub-Finsler Lie groups arise as solutions of a very natural class of optimal control problems \cites{control-system-GLn, control-systems-on-Lie-groups}. On the other hand, understanding geodesics in sub-Finsler Lie groups plays a key role in estimating the rate of convergence of nilpotent Lie groups to their asymptotic cones \cites{Pansu-croissanceboules, Breuillard-Enrico, libroenrico, bodart2025intermediate}. Asymptotic cones of finitely generated nilpotent groups equipped with word metrics are Carnot groups equipped with polyhedral norms \cite[Chapter 12]{libroenrico}. 
For the standard objects associated with sub-Finsler Lie groups, such as their polarizations and the norms on them, we refer to \cite{libroenrico} or \cref{sec:preliminaries}.
In this article, we do not assume any condition on the smoothness of the norm, nor do we require the norm to be strictly convex, as indeed in the polyhedral case. 
To study geodesics, i.e., isometric embeddings of intervals, we use the Pontryagin Maximum Principle, see \cref{thm:pmp}, a fundamental result in the theory of optimal control, which provides necessary conditions for curves to be geodesics. In this article, we focus on \textit{normal curves}, a particular class of curves that satisfy these necessary conditions. Normal curves are more manageable, since they satisfy a differential inclusion. We call the latter differential inclusion the \emph{normal inclusion}, see \eqref{eq:normal_equation}, since it generalizes the normal equation from sub-Riemannian geometry. 
The normal inclusion depends linearly on a covector in the dual of the Lie algebra of the group, which we refer to as a \textit{covector associated to the normal curve} (see \cref{def:normal+associated-covector}).

Normal curves are very well behaved when the norm is smooth and \textit{strongly convex}, in the standard sense of \cref{def:strictly-strongly-convex}. 
For example, if the norm is induced by a scalar product, the differential inclusion solved by normal curves becomes an analytic ordinary differential equation \cite[Theorem~7.3.3]{libroenrico}. Moreover, there exists an open dense set of the group with the property that every point of this set is connected to the unit element of the group by a unique and normal geodesic \cite[Section~11]{nuovo_libro_ABB}. 
More generally, if the norm is strongly convex, normal curves are of class $C^{1,1}$ (see \cref{prop:ham-reform-norm-eq}) and they are locally geodesic (up to linear reparametrization) \cite{Here-minimality-strongly-convex}. Moreover, for every point $p$ in the group and every covector $\lambda$ in the dual of the Lie algebra, there exists a unique normal curve starting at $p$ with associated covector $\lambda$
(see \cref{cor:uniqueness_sol_strongly_convex}).
Instead, in the case in which the norm is neither strictly convex nor smooth, e.g., when the norm is a \textit{polyhedral norm}, as in \cref{def:polytope-norm}, normal curves may become very wild: they may have corners, they might not be locally geodesics (see \cref{sec:non-locally-optimal-extremal-engel}),
and there may exist multiple normal curves starting at a point with the same associated covector. Nonetheless, normal curves always lift to normal curves via Lie homomorphisms that are submetries, see \cref{prop:lift-normal-curves}.

The first problem that we investigate in this article 
is the presence of branching of normal curves in sub-Finsler Lie groups. We say that in a sub-Finsler Lie group there is \emph{branching of normal curves} if there exist two different normal curves that coincide on a sub-interval of their domain, see \cref{def:branching}. 

When the norm is not strictly convex, branching may happen as a consequence of the non-uniqueness of the solution of the normal inclusion (e.g., this already happens in the normed vector space $(\R^2, \|\cdot\|_{L^1})$). However, when the norm is strongly convex, for every covector $\lambda$ there exists a unique normal curve with associated covector $\lambda$, and whether branching of normal curves may happen becomes unclear. In abelian groups, the presence of branching depends only on the convexity of the norm: as soon as the norm is strictly convex, all normal curves are lines, and branching of normal curves does not happen. Vice versa, if the norm is not strictly convex, then there exist multiple normal geodesics associated to the same covector, and there is branching of normal geodesics. Hence, one might think that the convexity of the norm captures the essence of the branching phenomenon.
However, in this article, we show that in the non-abelian setting, the smoothness of the norm also plays a big role.
On the one hand, we show that if the norm is strongly convex and the energy is continuously differentiable, then the branching of normal curves does not happen, see \cref{prop:Finsler-smooth}. On the other hand,
we show that as soon as the adjoint map is unbounded, there exists a (non-smooth) strongly convex norm for which the branching of normal geodesics happens. 

\begin{theorem}
\label{thm:if-condition-then-exists-branching}
    Let G be a Lie group with polarization $V$. Assume that there exists $X,Y\in V$ such that the set $\{ \Ad_{\exp(t Y)}X : t\in\R\}$ is unbounded. Then, there exists a strongly convex norm $\|\cdot\|$ on $V$ for which there is branching of normal geodesics in $G$ for the sub-Finsler metric given by $(V,\|\cdot\|)$.
\end{theorem}

The proof of \cref{thm:if-condition-then-exists-branching} is in \cref{sec:branching-proof-thm}. Moreover, in \cref{sec:example-2-dim-non-ab-group} we provide examples where we apply the theorem.

We emphasize that understanding whether there is branching of geodesics plays a role in the theory of mass transportation and hence in the study of synthetic curvature bounds in sub-Finsler Lie groups; see, for example, \cites{branching1, branching2}.

In the second part of this article, we consider norms for which the closed unit ball is a polyhedron (we refer to such norms as \emph{polyhedral norms}). We study whether normal curves exhibit some \textit{chattering}-type phenomenon. 
Namely, we investigate whether for every control $u$ of a normal curve, there exists $\delta>0$ such that the control $u$ takes value in a face on intervals of size smaller than $\delta$.
In this article, for a sub-Finsler Lie group $(G,V,\|\cdot\|)$
, we say that a measurable function $u:I\to V$ on an interval $I$ \emph{takes value in a face} if there exists $r>0$ and a face $A$ of the sphere $\partial B_{\|\cdot\|}(0,r)\subseteq V$ such that $u(t)\in A$ for almost every $t\in I$.
Here the \emph{size} of an interval $I=(a,b)$ is $|I|:=|a-b|$.

\begin{theorem}
\label{thm:main-polytope-norms-projection2}
    Let $G$ be a sub-Finsler Lie group with polyhedral norm $\|\cdot\|$. 
    \begin{itemize}
        \item[i.] For each normal curve $\gamma:[0,1]\to G$ there exists $\delta'>0$ such that for each interval $I\subseteq [0,1]$ with $|I|\leq \delta'$, the control of $\gamma|_{I}$ takes value in one face.  
        \item[ii.] More quantitatively, fixing a norm $N$ on the Lie algebra $\g$ of $G$ and a norm $N^*$ on the dual vector space $\g^*$, for $r>0$, set 
        \begin{equation}
        \label{eq:def-M}
            M(r):=\max\left\{N(\Ad_{g}\ad_XY) : 
        \begin{array}{c}
         g\in \bar B(1_G,r),\ X,Y\in\g , \\
   N(X)=N(Y)=1
        \end{array}
        \right\}.
        \end{equation}
        Then, there exists $\delta>0$ with the following property: for every $r>0$ and for each normal curve $\gamma:[0,1]\to G$ with $\gamma(0)=1_G$ and speed $r$, for each covector $\lambda$ associated to $\gamma$, and for every interval $I$ with $|I|\leq \frac{\delta}{N^*(\lambda)M\left(r\right)}$, the control of $\gamma|_I$ takes values in one face.
    \end{itemize}
\end{theorem}

The proof of \cref{thm:main-polytope-norms-projection2} is in \cref{subs689c584c}.

We remark that, although controls of normal curves take values in a face on intervals of small size, they may still not be continuous at an infinite number of points in some compact interval. Moreover, there exist sub-Finsler Carnot groups with polyhedral norms in which the unique (and normal) geodesic between two points has control that is non-continuous at an infinite number of points, see \cite{Lev2025}.

As a consequence of \cref{thm:main-polytope-norms-projection2}, we have that normal curves in Carnot groups have the property that their projection to the abelianization is locally normal and thus locally geodesic. Recall that, for every Carnot group $G$, its abelianization $G/[G,G]$ is a normed vector space and the projection $\pi:G\to G/[G,G]$ is a submetry, see \cite[Section 10.1.1]{libroenrico}.

\begin{corollary}
\label{cor:minimality-normal-curves}
Let $G$ be a Carnot group with polyhedral norm.
\begin{itemize}
    \item[i.] Normal curves in $G$ are locally length-minimizing.
    \item[ii.] More quantitatively, fixing a norm $N^*$ on the dual $\g^*$ of the Lie algebra of $G$, there exists $\alpha>0$ such that every normal curve $\gamma:[0,1]\to G$ of length $1$ and covector $\lambda\in\g^*$ is such that $\pi\gamma|_J$ is length-minimizing in $G/[G,G]$, for each interval $J\subseteq [0,1]$ satisfying $|J|\leq \frac{\alpha}{N^*(\lambda)}$.
\end{itemize}
\end{corollary}

The proof of \cref{cor:minimality-normal-curves} is in \cref{sec:carnot-polyhedral}.
We point out that the hypothesis in \cref{cor:minimality-normal-curves} of being in a Carnot group is crucial,
since we have an example of a sub-Finsler Lie group with polyhedral norm and with a normal curve that is not locally geodesic, see \cref{sec:non-locally-optimal-extremal-engel}.

When the polarization is the whole Lie algebra, we are dealing with Finsler Lie groups. In these groups, the norm of the covector associated to each normal curve is a multiple of its length. We use this fact to improve the quantitative estimate in \cref{thm:main-polytope-norms-projection2}. 

\begin{corollary}
\label{cor:finsler-compact}
    Let $G$ be a Finsler Lie group with polyhedral norm. Then there exists $L>0$ such that all controls of normal curves of length smaller than $L$ take values in a face. 
\end{corollary}

The proof of \cref{cor:finsler-compact} is in \cref{subs689c584c}. This last result does not extend to the sub-Finsler case, e.g., in the Heisenberg group as in \cite{l1-heisenberg-reference}.

\subsection{The key idea}
Let $(G,V,\|\cdot\|)$ be a sub-Finsler Lie group, where $V$ is a polarization and $\|\cdot\|$ is a norm on $V$. Both \cref{thm:if-condition-then-exists-branching} and \cref{thm:main-polytope-norms-projection2} rely on the following simple but fundamental observation: in the normal inclusion
 \begin{equation*}       \lambda\Ad_{\gamma(t)}|_V\in \partial_{u(t)}\E,\quad \text{for a.e. } t\in I, 
\end{equation*}
the left-hand side is an absolutely continuous curve that takes value in a sphere in $V^*$ (see \cref{prop:subdifferentials-of-energy}). Hence, if the curve $\lambda\Ad_{\gamma}|_V$ is inside an open set of such a sphere at some time $t$, it will stay inside the same open set in an open neighborhood of $t$. 

The idea behind \cref{thm:if-condition-then-exists-branching} is to define a norm on $V$ that has a corner at $Y$. To be more precise, we require the subdifferential of the energy at $Y$ to have non-empty interior $U$ in the dual sphere of radius $\|Y\|$.
 Then, if a curve $\gamma$ solves the normal inclusion with covector $\lambda\in U$, we have $\lambda\Ad_{\gamma}|_V\in U $ on some open interval, thus $\gamma$ will coincide with the one-parameter subgroup $t\mapsto \exp(tY)$ on this interval. The condition on the unboundedness of the adjoint will imply the existence of a time $t$ for which $\lambda\Ad_{\gamma(t)}|_V\notin U $. Thus, at time $t$, the curve $\gamma$ does not coincide anymore with a one-parameter subgroup. Using this strategy, we find a family of normal curves that coincide with $t\mapsto \exp(tY)$ on open intervals of different sizes, proving the presence of branching.

When instead we deal with polyhedral norms, we find an open covering of the dual sphere with open sets that have the following property: if the curve $\lambda\Ad_{\gamma}|_V$ is valued into one of these open sets on some interval, then on the same interval the control of $\gamma$ stays on a face of the polyhedron. We then find some $\delta'>0$ such that on all intervals of size smaller than $\delta'$, the curve $\lambda\Ad_{\gamma}|_V$ is valued into one of such open sets. To do so, we give an estimate of the dual norm of the derivative of $\lambda\Ad_{\gamma}|_V$.

\subsection{Outline of the paper}
In \cref{sec:preliminaries} we review some notions of convex analysis that will be useful in the paper. We define strongly convex functions and norms, and we present Fenchel conjugates and their properties. We then present several properties of the subdifferential of the energy function. We finish the section with a quick overview of polyhedral norms and some of their properties.
In \cref{sec:sub-FinslerLiegroups} we introduce the terminology of sub-Finsler Lie groups. We present the Pontryagin Maximum Principle, and we define normal curves. We prove the existence of solutions of the normal inclusion and that the solution is unique if the norm is strongly convex (see \cref{prop:solution-normal-exists} and \cref{cor:uniqueness_sol_strongly_convex}). At the end of the section, we briefly recall the definition of Carnot groups.
In \cref{sec:branching} we prove that for $C^1$ strongly convex norms there is no branching of geodesics (see \cref{prop:Finsler-smooth}). We then give a sufficient condition for branching of normal curves to happen in sub-Finsler Lie groups (see \cref{them:main-branching-of-normal}). \cref{thm:if-condition-then-exists-branching} will follow from the latter sufficient condition. At the end of section 4, we present two examples. In the first one, we consider the affine group of the line with an explicit sub-Finsler norm, to help the reader familiarize with the ideas in the proof of \cref{them:main-branching-of-normal}. Instead, the second example is in $SO(3)$ and shows that the condition on the adjoint is needed in \cref{them:main-branching-of-normal}.
In \cref{sec:polytope}, we study sub-Finsler Carnot groups with polyhedral norms and we prove \cref{thm:main-polytope-norms-projection2}, \cref{cor:minimality-normal-curves}, and \cref{cor:finsler-compact}. Finally, we give an example of a normal extremal that is not locally optimal in the Heisenberg group equipped with some Finsler norm.

\subsection{Acknowledgements}
We would like to thank Lev Lokutsievskiy and Michele Motta for the insightful discussions on the topic of this paper. We would like to thank Lev~Lokutsievskiy for inspiring the example presented in \cref{sec:non-locally-optimal-extremal-engel}.

\section{Preliminaries in convex analysis}
\label{sec:preliminaries}

In this section, we introduce some notions and results of convex geometry. We mainly follow \cites{Clarke-book, rockafellar2015convex}. 

\subsection{Elements from convex analysis}
In this section, we introduce a few elements from convex analysis.
First, we recall the definitions of strict and strong convexity.
Next, we recall the definitions of the Fenchel dual and of the subdifferential of a function.
We conclude with a few properties linking these three notions.
All the material in this section is classical.

\begin{definition}\label{def:strictly-strongly-convex}
    Let $V$ be a vector space. 
    A function $f:V\to \R$ is \emph{convex} if 
    \begin{equation}
        f(tv+(1-t)u) \le tf(v)+(1-t)f(u), \quad \forall u,v\in V,\ \forall t\in(0,1). 
    \end{equation}
    A function $f:V\to \R$ is \emph{strictly convex} if 
    \begin{equation}
    \label{eq:strictly-convex-function}
        f(tv+(1-t)u)< tf(v)+(1-t)f(u), \quad \forall u,v\in V,\text{ with } u\neq v ,\forall t\in(0,1). 
    \end{equation}
    A function $f:V\to \R$ is \emph{strongly convex} if there exists a euclidean norm $N:V\to \R$ such that for all $u,v\in V$ and all $t\in[0,1]$ there holds
    \begin{equation}
    \label{eq:strongly-convex-function}
        f(tv+(1-t)u)\leq tf(v)+(1-t)f(u)-t(1-t)\left(N(u-v)\right)^2. 
    \end{equation}
\end{definition}

In this section, we denote by $V$ a normed vector space, for example, the Hilbert space $L^2$ or a finite-dimensional normed vector space. For such a $V$, we denote by $V^*$ the space of linear continuous functions from $V$ to $\R$ and by $\langle\cdot,\cdot\rangle:V\times V^*\to \R$ the duality pairing. Recall that we have a standard canonical inclusion $V\subseteq V^{**}$, which is an equality when $V$ is finite-dimensional. In this article, we use notation from \cite{Clarke-book}, and for example we consider $\R_{\infty} := \R\cup \{+\infty\}$.

\begin{definition}[{Fenchel conjugates, \cite[$\S$4.2]{Clarke-book}}]
\label{def:fenchel-conjugate}
Let $V$ be a normed vector space and $f:V \to \R_\infty$ such that $f\not \equiv +\infty$. The \emph{Fenchel conjugate} of $f$ is the map $f^*: V^*\to \R_{\infty}$ defined by setting
\begin{equation*}
    f^*(\eta):=\sup_{v\in V} \left(\langle\eta,v\rangle-f(v)\right),\quad  \forall \eta\in V^*.
\end{equation*}
\end{definition}

The Fenchel conjugate is well behaved when the function is continuous and convex, as recalled in the next propositions.

\begin{proposition}\emph{(Fenchel-Moreau, \cite[Theorem 4.21]{Clarke-book}). }
    Let $V$ be a finite-dimensional vector space and $f:V\to \R_\infty$ such that $f\not \equiv +\infty$. We use the canonical identification $V^{**}=V$. Then, the function $f$ is convex lower-semi-continuous
    if and only if $f^*\not \equiv +\infty$ and $f=f^{**}$. 
\end{proposition}

There is a characterization of strict and strong convexity in terms of Fenchel conjugates.
If $f:V\to\R$ is differentiable at $v\in V$, we denote by $\df_vf$ the derivative of $f$ at $v$.
Since $\df_vf$ is a linear map $V\to\R$, we have $\df_vf\in V^*$.

\begin{proposition}[{Properties of Fenchel conjugates  \cite[Chapter X, Section~4]{referenza-C11-Fenchel-conjugate}}]\label{prop:properties-fenchel-conj}
 Let $V$ be a finite-dimensional vector space and $f:V\to \R$ be a convex continuous function. Then, $f^*\in C^1(V^*)$ if and only if $f$ is strictly convex. Moreover, $f^*\in C^{1,1}(V^*)$ if and only if $f$ is strongly convex. When $f$ is strictly convex, there holds
    \begin{equation}
\label{eq:differential-fenchel-conjugate}
        \langle \df_\eta f^*, \xi\rangle=\langle \xi, \operatorname{argmax}_{v\in V}\left(\langle\eta,v\rangle-f(v)\right) \rangle, \quad \forall\eta,\xi\in V^* . 
    \end{equation}
\end{proposition}

The last notion that we take from convex analysis is the subdifferential of a function.
This notion of derivative is well-suited for convex functions.

\begin{definition}[{Subdifferential, \cite[$\S$4.1]{Clarke-book}}]
Let $V$ be a normed vector space, $f:V\to \mathbb{R}_\infty$, and $u\in V$ such that $f(u)\neq +\infty$. 
An element $\eta \in V^*$ is a \emph{subgradient} of $f$ at $u$ if
\begin{equation}
    \eta(v-u)\leq f(v)-f(u), \quad \forall v\in V.
\end{equation}
The set of subgradients at $v$ is denoted by    $\partial_uf\subseteq V^*$ and called the \emph{subdifferential} of $f$ at $u$.
\end{definition}
Subdifferentials of  convex functions at a point are   nonempty convex compact sets \cite[Corollary~4.7]{Clarke-book}. 
By \cite[Corollary~4.4]{Clarke-book}, if $f:V\to \R$ is convex, then \begin{equation}
\label{eq:subdifferential-gateaux-diff}
F \text{  G\^ateaux differentiable at  } u\in V \implies 
\partial_uf = \{\df_uf\}.\end{equation}
By \cite[Theorem 4.10]{Clarke-book}, if $f,g:V\to \mathbb{R}$ are continuous and convex, then   
\begin{equation}
\label{eq:sum-of-subdifferentials}
    \partial_xf + \partial_xg = \partial_x(f+g), \qquad \forall x\in V.
\end{equation}  

\begin{remark}[Subdifferential inversion, {\cite[Exercise 4.27]{Clarke-book}}]
\label{remark:duality-subdifferentials}
 Let $V$ be a normed vector space and $f:V\to\R_\infty$
  convex and lower semicontinuous with $f\not\equiv \infty$.
 Then, for every $u\in V$ and $\eta\in V^*$ with $f(u), f^*(\eta)\neq +\infty$, 
\begin{equation}\label{eq68a5df0a}
 \eta\in \partial_u f
 \quad\Leftrightarrow\quad
 f(u) + f^*(\eta) = \langle \eta,u \rangle
 \quad\Leftrightarrow\quad
 u\in\partial_\eta f^* .
\end{equation}
\end{remark}

\subsection{Strict and strongly convex norms}

We apply the above notions from convex analysis to the energy of norms.
\begin{definition}\label{def68a5d473}
    Let $\|\cdot\|$ be a norm on a vector space $V$. The \emph{energy} associated to $\|\cdot\|$ is the function $\E:V\to[0,+\infty)$, $\E(v) := \frac{1}{2}\|v\|^2$ for $v\in V$. 
    The norm $\|\cdot\|$ is \textit{strictly convex} (resp. \emph{strongly convex}) if $\E$ is a strictly convex (resp. strongly convex) function.
\end{definition}

\begin{remark}
 Definition~\ref{def68a5d473} overloads the expressions ``strictly convex'' and ``strongly convex''.
 However, notice that a norm is never ``strictly convex'' as a convex function in the sense of Definition~\ref{def:strictly-strongly-convex}.
 Indeed, if $v\in V$ is a nonzero vector and $t\in(0,1)$, then $\|v\| = \|(t+1-t)v\| = t\|v\| + (1-t)\|v\|$.
 So, when referred to norms, strict and strong convexity are in the sense of Definition~\ref{def68a5d473}.
 Such overloading comes from a by-now-standard lexicon.
\end{remark}

\begin{remark}
 Let $\|\cdot\|$  be a norm on a finite-dimensional vector space $V$.
 Then $\|\cdot\|$ is strictly convex if and only if the boundary of the unit ball does not contain segments. 
 If $\|\cdot\|$ is smooth outside the origin, then $\|\cdot\|$ is strongly convex if and only if the boundary of the unit ball does not contain points where the curvature vanishes.
\end{remark}

Let $V$ be a vector space and $\|\cdot\|$ a norm on $V$. On the space $V^*$ we define the \emph{dual norm} $\|\cdot\|_*$ by setting
\begin{equation}
\label{eq:def-dual-norm}
    \|\eta\|_*:=\sup\{\langle\eta,v\rangle : v\in V, \ \|v\|=1\}, \quad \forall \eta\in V^*.
\end{equation}
For $r\in[0,+\infty)$, we denote with $S(r)$ and $S^*(r)$ the sphere of radius $r$ in $V$ and $V^*$, respectively. When $r=1$, we write $S$ or $S^*$.

%

The following proposition describes the subdifferential of the energy at a point.
This is a simple exercise; see some details in \cite[Lemma 2.19]{Hakavuori-2020-step2_geodesics}.  

\begin{proposition}[{\cite[Exercise~4.2]{Clarke-book}}]
\label{prop:subdifferentials-of-energy}
    Let $(V,\|\cdot \|)$ be a normed vector space
    with energy $\E:=\frac{1}{2}\|\cdot\|^2$. 
    The Fenchel conjugate of $\E$ is $\E^*=\frac{1}{2}\left(\|\cdot\|_*\right)^2$.
    Moreover, for every $u\in V$ there holds
    \begin{equation}\label{eq:subdiff-of-energy}
    \begin{aligned}
  \partial_u\E &= \{\eta\in V^* : \|\eta\|_*=\|u\|,\ \langle \eta ,u \rangle=\|u\|^2\} \\
  &= \{\eta\in V^* : \|\eta\|_*=\|u\|,\ \langle \eta ,u \rangle=\|\eta\|_*\|u\|\} .
    \end{aligned}
    \end{equation}
    Symmetrically, for every $\eta\in V^*$, 
    \begin{equation}\label{eq:def_tan}
    \begin{aligned}
     \partial_\eta \E^* 
     &= \{v\in V : \|v\|=\|\eta\|_*,\ \langle \eta,v \rangle = \|\eta\|_*^2 \} \\
     &= \{v\in V : \|v\|=\|\eta\|_*,\ \langle \eta,v \rangle = \|\eta\|_*\|v\| \}.
 \end{aligned}
 \end{equation}
\end{proposition}
\section{Sub-Finsler Lie groups}
\label{sec:sub-FinslerLiegroups}

In this section, we state the Pontryagin maximum principle for sub-Finsler Lie groups, and we study normal curves. We make use of the notions of convex analysis presented in \cref{sec:preliminaries}. We shall use standard terminology from Lie group theory, as in \cite{libroenrico}.

\subsection{Normal curves in sub-Finsler Lie groups}
On a Lie group,
we denote with $L,R$ and $\Ad$ the left-translation, the right translation, and the adjoint map.
A \emph{polarization} of a Lie algebra $\g$ is a linear subspace of $\g$.
By extension, a \emph{polarization} of a Lie group is a polarization of its Lie algebra.
A polarization $V$ of a Lie algebra $\g$ is \emph{bracket-generating} if the Lie algebra generated by $V$ is the whole $\g$.
Let $G$ be a Lie group with Lie algebra $\g$ and polarization $V\subseteq \g$.
For every interval $I\subseteq \R$,
an absolutely continuous curve $\gamma:I\to G$ is \emph{horizontal} if 
\begin{equation*}
    \dot \gamma\in \df L_{\gamma(t)}V ,\quad\text{for a.e. } t\in[0,1].
\end{equation*}
In this case, we say that $u := \df L_ {\gamma}^{-1}\dot \gamma\in L^1([0,1],V)$ is the \emph{control} of $\gamma$.

For each norm $\|\cdot\|$ on $V$, we consider the \emph{Carnot-Carathéodory distance} $d:G\times G\to [0,+\infty]$ defined by, for $g,h\in G$, 
\begin{equation}
\label{eq:distance}
    d(g,h):=\inf\left\{ \int_0^1 \|\df L_{\gamma}^{-1}\dot \gamma(t)\|\df t  \ \ \Bigg| \ \  \begin{gathered}
       \gamma:[0,1]\to G \text{ horizontal,}\\
       \gamma(0)=g,\gamma(1)=h
    \end{gathered}  \right\}.
\end{equation}

By Chow--Rashewskii's Theorem, if $V$ is bracket generating, then $d$ is a finite-valued distance that induces the manifold topology, see \cite[Section~7.1.4]{libroenrico}.
\begin{definition}
    A \textit{sub-Finsler Lie group} is a triple $(G,V,\|\cdot\|)$, where $G$ is a Lie group,
    $V$ is a bracket-generating polarization of $G$,
    and $\|\cdot\|$ is a norm on $V$.
    When $V$ is full-dimensional, i.e., it equals the Lie algebra of $G$, we say that $(G,V,\|\cdot\|)$ is a \emph{Finsler} Lie group.
\end{definition}

Curves realizing the infimum in \eqref{eq:distance} are called \emph{length-minimizers}. 
When re\-para\-me\-trized with constant speed, length-minimizers are also \emph{energy minimizers}, that is, they also realize the infimum
\begin{equation*}
    \inf\left\{ \frac{1}{2}\int_0^1 \|\df L_{\gamma}^{-1}\dot \gamma(t)\|^2\df t  \ \ \Bigg| \ \  \begin{gathered}
       \gamma:[0,1]\to G \text{ horizontal,}\\
       \gamma(0)=g,\gamma(1)=h
    \end{gathered}  \right\}.
\end{equation*}

Up to affine reparametrization, energy minimizers are \emph{geodesics}, that is, isometric embedding of intervals.
The Pontryagin Maximum Principle gives necessary conditions for curves to be energy minimizers. We rephrase the Pontryagin Maximum Principle using the terminology introduced in \cref{sec:preliminaries}.

\begin{theorem}[Pontryagin Maximum Principle]
\label{thm:pmp}
 Let $(G,V,\|\cdot\|)$ be a sub-Finsler Lie group with
 Lie algebra $\g$ and
  energy $V\ni v \mapsto \E(v) := \frac{\|v\|^2}{2}\in\R$.
 If a horizontal curve $\gamma:[0,1]\to G$ with control $u:[0,1]\to V$ is an energy minimizer, then there exists $\lambda\in \g^*$ such that either
 \begin{equation}
         \lambda\Ad_{\gamma(t)}|_V\in \partial_{u(t)}\E, \ \ \text{ for a.e. } t\in[0,1],
         \label{e normale}
 \end{equation}
 or $\lambda\neq 0$ and
 \begin{equation}
     \lambda\left(\Ad_{\gamma(t)}X\right)
     =0, \ \ \forall t\in[0,1] , \ \forall  X\in V.
     \label{eq:abnormale}
 \end{equation}
\end{theorem}

Although one can find several proofs of the Pontryagin Maximum Principle in the literature, we describe here a simplified version that
uses a theorem by Clarke \cite[Theorem 10.47]{Clarke-book}
and it is adapted to sub-Finsler Lie groups; see also \cite[Proposition 2.14]{LDP-escape}.

\begin{proof}[Proof of \cref{thm:pmp}]
 Let $\End:L^2([0,1],V)\to G$ be the map  $\End(u):=\gamma_u(1)$, where $\gamma_u:[0,1]\to G$ is the unique curve with control $u$ starting from the identity element of the group $G$.
 Let $\E_{L^2}:L^2([0,1],V)\to \R$ be the map defined by $\E_{L^2}(u):=\frac12\int_0^1\|u(t)\|^2\df t = \int_0^1\E(u(t)) \df t$, for all $u\in L^2([0,1],V)$.
 
 Let $\gamma:[0,1]\to G$ be a horizontal curve with control $u$ that is energy minimizing. Up to a left-translation, we assume $\gamma(0)$ to be the identity element of $G$. Since $\gamma$ is energy minimizing, we have $u\in L^2([0,1],V)$.
 By Clarke's version of the Lagrange Multiplier Theorem \cite[Theorem 10.47]{Clarke-book}, we get that there exists $\bar\lambda\in T^*_{\gamma(1)}G$ and $\nu\in \{0,1\}$ such that $(\bar \lambda,\nu)\neq 0$ and
 \begin{equation}\label{eq689c5e37}
     0\in \partial_u\left(\bar\lambda\End+\nu\E_{L^2}\right).
 \end{equation}
 Since $\End$ is a smooth function (as proved in~\cite{Libro-Rifford}), using \eqref{eq:sum-of-subdifferentials} and 
 \eqref{eq:subdifferential-gateaux-diff} we get that~\eqref{eq689c5e37} is equivalent to
 \begin{equation*}
   -\bar\lambda \df_u\End\in \nu\partial_u\E_{L^2}.
 \end{equation*}
 Using the explicit formula for the differential of the end-point map in \cite[Proposition~7.2.1]{libroenrico}, and setting $\lambda:=-\bar\lambda \df R_{\gamma(1)}^{-1}$, we obtain for each $v\in L^2([0,1],V)$ the inequality
 \begin{equation}\label{eq:prima-formula-proof-pmp}
     \int_0^1\lambda(\Ad_{\gamma(t)} v(t))-\int_0^1\lambda(\Ad_{\gamma(t)} u(t))\df t 
     \leq \nu \int_0^1 (\E(v(t))-\E(u(t)))\df t. 
 \end{equation}
 
 Fix $Y\in V$. Let $\bar{t}$ be a Lebesgue point of $u$. 
 For each $\delta>0$, we apply~\eqref{eq:prima-formula-proof-pmp} to the test function 
 $v:= u-(Y-u(\bar t))\chi_{[\bar{t},\bar t+\delta]}$, where $\chi_{[\bar{t},\bar t+\delta]}$ is the charachteristic function of the interval $[\bar{t},\bar t+\delta]$.
 Therefore, we obtain that, for every $\delta>0$,
 %
 \begin{equation}\label{eq689c615c}
     \frac{1}{\delta}\int_{\bar{t}}^{\bar{t}+\delta}\lambda(\Ad_{\gamma(t)} (Y-u(\bar{t}))\df t
     \leq \frac{\nu}{\delta} \int_{\bar{t}}^{\bar{t}+\delta} (\E(u(t)-Y-u(\bar t))-\E(u(t)))\df t.
 \end{equation}
 Since $\E$ is Lipschitz, $\bar{t}$ is also a Lebesgue point of $t\mapsto \E(u(t))$ and of $t\mapsto\E(u(t)-Y-u(\bar t))$. 
 Hence, passing to the limit $\delta\to0^+$ in~\eqref{eq689c615c}, we conclude
 \begin{equation*}
    \lambda(\Ad_{\gamma(\bar t)} (Y-u(\bar{t})) \leq  \nu \left(\E(Y)-\E(u(\bar t))\right).
 \end{equation*}
 Since $Y\in V$ was arbitrary, we have proved
 \begin{equation}
 \label{eq:proof-pmp-seconda}
    \lambda\Ad_{\gamma(\bar t)}|_V \in \nu \partial_{u(\bar t) }\E.
 \end{equation}
 If $\nu=1$, the inclusion \eqref{eq:proof-pmp-seconda} rewrites as \eqref{e normale}; when $\nu=0$ it rewrites as \eqref{eq:abnormale}.
\end{proof}

\begin{definition}
\label{def:normal+associated-covector}
    Let $(G,V,\|\cdot\|)$ be a sub-Finsler Lie group with energy $\E:V\to\R$, and let $I\subseteq \R$ be an interval. 
    A horizontal curve $\gamma:I\to G$ with control $u\in L^2(I,V)$ is a \emph{normal curve} if there exists $\lambda\in \g^*$ such that
    \begin{equation}\label{eq:normal_equation}
        \lambda\Ad_{\gamma(t)}|_V\in \partial_{u(t)}\E,\quad \text{for a.e. } t\in I.
 \end{equation}
 Each covector $\lambda\in \g^*$ for which the inclusion \eqref{eq:normal_equation} holds is called a (normal) \emph{covector associated to} $\gamma$.
\end{definition}

\begin{remark}[Normal inclusion as differential inclusion]
\label{rem:normal-differential-inclusion}
    After the subdifferential inversion \eqref{eq68a5df0a}, 
    we can rewrite \eqref{eq:normal_equation} as a differential inclusion
    \begin{equation}
        \label{eq:normal-differential-inclusion}
        u(t)\in \partial_{\lambda\Ad_{\gamma(t)}|_V} \E^* ,\quad \text{for a.e. }t\in I,
    \end{equation}
    where $\E^*:V^*\to\R$ is the energy of the dual norm $\|\cdot\|_*$ on $V^*$,
    which is the Fenchel conjugate of $\E$ by \cref{prop:subdifferentials-of-energy}.
\end{remark}

\begin{proposition}[Existence of solutions to the normal inclusion]
 \label{prop:solution-normal-exists}
 Let $G$ be a sub-Finsler Lie group with Lie algebra $\g$.
 Then for all $\lambda\in \g^*$ and all $g\in G$ there exists a normal curve starting at $g$ with associated covector $\lambda$.
\end{proposition}
\begin{proof}
 Let $V$ be the polarization of $G$, $\|\cdot\|$ the norm, and $\E$ its energy.
 Being $\E^*$ convex, by \cite[Corollary 4.7]{Clarke-book}, the set $\partial_{\lambda\Ad_{x}|_V} \E^*$ is non-empty, convex, and compact. 
 Moreover, by \cite[Proposition 4.14]{Clarke-book}, the set-valued function $x\mapsto \partial_{\lambda\Ad_{x}|_V} \E^*$  is upper semicontinuous. 
 Hence, by Filippov's Theorem \cite[Theorem 1, Pag 77]{Filippov}, we have that for all $\lambda\in \g^*$ and all $g\in G$ there exists a solution $\gamma$ of the differential inclusion \eqref{eq:normal-differential-inclusion} with $\gamma(0)=g$.
\end{proof}

Left translations and affine reparametrization of normal curves are normal curves; see~\cite[Remark 2.19]{LDP-escape}. 

\begin{remark}\label{rem68ad82a9}
 When the norm is strongly convex, it is known that normal curves are locally energy minimizing; see~\cite[Section 17]{agrachev2013control} or~\cite{Here-minimality-strongly-convex}.  
 However, there are examples of normal curves that are not locally length-minimizing, see our example in \cref{sec:non-locally-optimal-extremal-engel}.
\end{remark}

It is known that normal curves are parametrized by multiple of arc-length:

\begin{proposition}\emph{(\cite[Proposition 2.23]{LDP-escape}).}\label{prop:normal-PBAL}
    Normal curves in sub-Finsler Lie groups are parametrized with constant speed equal to $\|\lambda|_V\|_*$, where $\lambda$ is any associated covector and $V$ is the polarization.
\end{proposition}

\subsection{Normal curves for strongly convex norms}
When the norm in a sub-Finsler Lie group is strongly convex, 
the differential inclusion~\eqref{eq:normal-differential-inclusion} becomes a differential equation that has uniqueness of solutions.


\begin{proposition}\label{prop:ham-reform-norm-eq}
 Let $G$ be a sub-Finsler Lie group equipped with a strictly convex norm with energy $\E:V\to \R$. 
 Let $\gamma:[0,1]\to G$ be a horizontal curve with control $u$. 
 Then $\gamma$ is normal if and only if there exists a covector $\lambda\in \g^*$ such that    
 \begin{equation}
 \label{eq:normal-equation-ham-version}
     u(t)  =\df_{\lambda\Ad_{\gamma(t)}|_V} \E^*,\quad  \text{ for a.e. } t\in [0,1].  
 \end{equation}
 where $\E^*:V^*\to\R$ is the Fenchel conjugate of the energy $\E$,
 or, equivalently, $\E^*$ is the energy of the dual norm.
\end{proposition}
\begin{proof}
 The fact that the Fenchel conjugate of the energy $\E$ is equal to the energy of the dual norm follows from \cref{prop:subdifferentials-of-energy}. Since by \cref{prop:properties-fenchel-conj} the Fenchel conjugate $\E^*$ of the strictly convex function $E$ is $C^1$, by \eqref{eq:subdifferential-gateaux-diff} the normal inclusion \eqref{eq:normal-differential-inclusion} becomes the differential inclusion \eqref{eq:normal-equation-ham-version}.
\end{proof}

As an immediate consequence, we get the following uniqueness result.

\begin{proposition}[Uniqueness of solution of the normal equation for strongly convex norms]
\label{cor:uniqueness_sol_strongly_convex}
 Let $G$ be a sub-Finsler Lie group endowed with a strongly convex norm.
 Denote by $\g$ the Lie algebra of $G$. 
 For all $p\in G$ and $\lambda\in \g^*$, there exists a unique normal curve $\gamma:[0,1]\to G$ with associated covector $\lambda$ satisfying $\gamma(0)=p$. Moreover, we have $\gamma\in C^{1,1}([0,1],G)$. 
\end{proposition}
\begin{proof}
 By \cref{prop:ham-reform-norm-eq}, every normal curve is a solution of the ODE \cref{eq:normal-equation-ham-version}, 
 whose right-hand side is Lipschitz by \cref{prop:properties-fenchel-conj}.
 By standard ODE results, \cref{eq:normal-equation-ham-version} has unique solutions.
\end{proof}

\subsection{Lifts of normal curves via submetries}

In this section, we show that whenever one has a Lie group epimorphism $G\to H$, 
every sub-Finsler structure on $G$ can be pushed forward to a sub-Finsler structure on $H$ such that the epimorphism is a metric submetry. 
Moreover, normal curves in $H$ are projections of normal curves in $G$.
In what follows, we revise the notion of  submetry between metric spaces, 
then we explain a couple of observations for submetries between normed vector spaces,
and finally, we prove \cref{prop:lift-normal-curves} for the lift of normal curves between sub-Finsler Lie groups.
Recall that length-minimizing curves lift to length-minimizing curves via submetries in abstract metric spaces; see \cite[Corollary 3.1.25]{libroenrico}.

%


\begin{definition}
\label{def:submetry}
    A map $f:G\to H$ between metric spaces is a \emph{submetry} if 
    \begin{equation*}
        f(\bar{B}(p,r))=\bar{B}(f(p),r),\quad \forall p\in G,\forall r>0,
    \end{equation*}
    where with $\bar B(p,r)$ we denote the closed ball of radius $r$ centered at $p$.
\end{definition}

For a linear map $A:V\to W$ between normed vector spaces, we defined the \textit{transpose} $A^*:W^*\to V^*$ of $A$ by setting $A^*\eta = \eta\circ A$ for all $\eta\in W^*$.
{The following proposition is a characterization of submetries between normed vector spaces.}

\begin{proposition}
\label{prop:dual-submetry-is-iso-embedding}
        Let $V$ and $W$ be normed vector spaces and $A:V\to W$ a linear bounded operator.
        If $A$ is a submetry, then $A^*:W^*\to V^*$ is an isometric embedding.
        
        Suppose in addition that $V$ is a reflexive Banach space.
        Then, $A$ is a submetry if and only if $A^*:W^*\to V^*$ is an isometric embedding.
\end{proposition}

\begin{proof}
        Assume that $A$ is a submetry. Let $\beta \in W^*$. Then
        \begin{align*}
                \|A^*\beta\|_{V^*}
                &= \sup\{ \langle A^*\beta , v \rangle : v\in V,\ \|v\|_V\leq 1\} \\
                &= \sup\{ \langle \beta , Av \rangle : v\in V,\ \|v\|_V\leq 1\} \\
                &= \sup\{ \langle \beta , w \rangle : w\in W,\ \|w\|_W\leq 1\} \\
                &= \|\beta\|_{W^*},
        \end{align*}
        where in the third equality we used the fact that $A(\bar B(0,1)) = \bar B(0,1)$.
        
        Vice versa, assume that $V$ is a reflexive Banach space and that $A^*:W^*\to V^*$ is an isometric embedding. Let $v \in V$. Then
        \begin{align*}
                \|Av\|_W
                &\stackrel{\eqref{eq:def-dual-norm}}{=} \sup\{\langle \beta, Av \rangle : \|\beta\|_{W^*} \leq 1\} \\
                &= \sup\{\langle A^*\beta, v \rangle : \|\beta\|_{W^*} \leq 1\} \\
                &\le \sup\{\langle \alpha, v \rangle : \|\alpha\|_{V^*} \leq 1\} \\
                &= \|v\|_V.
        \end{align*}
        Therefore, $A(\bar{B}_V(0,1)) \subset \bar{B}_W(0,1)$.
       
       Regarding the opposite inclusion, we prove it in two steps.
       First, we claim that $\bar{B}_W(0,1) \subset\mathrm{closure}(A(\bar B_V(0,1)))$.
       Indeed, by the Hahn–Banach Theorem~\cite{Rudin}, if $w\in \bar{B}_W(0,1) \setminus \mathrm{closure}(A(\bar B_V(0,1)))$, then 
%
        there exists $\beta \in W^*$ such that $\mathrm{closure}(A(\bar B_V(0,1))) \subset \{\beta < 1\}$ and $\langle \beta, w \rangle > 1$.
        But then
        \begin{align*}
                \|A^*\beta\|_{V^*}
                &\stackrel{\eqref{eq:def-dual-norm}}{=} \sup\{\langle A^*\beta, v \rangle : v \in V,\ \|v\|_V \leq 1\} \\
                &= \sup\{\langle \beta, Av \rangle : v \in V,\ \|v\|_V \leq 1\} \\
                &< 1
                < \langle \beta, w \rangle \\
                &\leq \sup\{\langle \beta, z \rangle : z \in W,\ \|z\|_W \leq 1\}
                \stackrel{\eqref{eq:def-dual-norm}}{=} \|\beta\|_{W^*},
        \end{align*}
        which contradicts the assumption that $A^*$ is an isometric embedding.
        We have thus shown that $\bar{B}_W(0,1) \subset\mathrm{closure}(A(\bar B_V(0,1)))$.
       
       Second, we claim that $\bar{B}_W(0,1) \subset A(\bar B_V(0,1))$.
       Here we use the reflexivity of $V$.
       Indeed, for each $w\in\bar{B}_W(0,1)$, by the previous claim there exists $\{v_n\}_{n\in\N}\subset \bar B_V(0,1)$ such that $\lim_{n\to\infty} Av_n=w$.
       Since $V$ is reflexive, by the Banach--Alaouglu Theorem~\cite{Rudin}, up to passing to a subsequence, there exists $v_\infty\in \bar B_V(0,1)$ with $v_n\rightharpoonup v_\infty$ as $n\to\infty$.
       Therefore, for every $\beta\in W^*$,
       \begin{align*}
        \langle \beta,w \rangle
&= \lim_{n\to\infty} \langle \beta,Av_n \rangle
= \lim_{n\to\infty} \langle A^*\beta,v_n \rangle
= \langle A^*\beta,v_\infty \rangle
= \langle \beta,Av_\infty \rangle .
       \end{align*}
       We conclude that $Av_\infty=w$ and thus $w\in A(\bar B_V(0,1))$.
\end{proof}

As a consequence of \cref{prop:dual-submetry-is-iso-embedding}, we get that  the transpose of a submetry maps subgradients of the energy to subgradients of the energy:

\begin{corollary}\label{cor:lift-subdifferential}
    Let $A:V\to W$ be a linear submetry between normed vector spaces.
    We denote by $\|\cdot\|$ both norms of $V$ and $W$, and by $\E$ both energies of $V$ and $W$. 
 Let $w\in W$ and $\eta \in \partial_w\E$. 
 Then, for all $v\in A^{-1}(w)$ satisfying $\|v\|=\|w\|$, we have $A^*\eta \in \partial_v\E$.
\end{corollary}
\begin{proof}
    Fix $v\in A^{-1}(w)$ such that $\|v\|=\|w\|$. 
    Hence, by \cref{prop:subdifferentials-of-energy} and \cref{prop:dual-submetry-is-iso-embedding} we have 
    \begin{equation}
       \label{eq:norma-duale-pull-back-via-A}\|A^*\eta\|_*=\|\eta\|_*\stackrel{\eqref{eq:subdiff-of-energy}}{=}\|w\|=\|v\|. 
    \end{equation}
Moreover, we have 
\begin{eqnarray}
    \label{eq:verifica-norma-realizzata-da-v-in-prova-pullback}
    \|\langle A^*\eta,v\rangle\|=\|\langle\eta,A(v)\rangle\|=\|\langle\eta,w\rangle\|
\stackrel{\eqref{eq:subdiff-of-energy}}{=}\|w\|^2=\|v\|^2.
\end{eqnarray}
Thus, by Equations \eqref{eq:norma-duale-pull-back-via-A}-\eqref{eq:verifica-norma-realizzata-da-v-in-prova-pullback} and \cref{prop:subdifferentials-of-energy} we conclude $A^*\eta\in \partial_v\E$.   
\end{proof}

We now focus on sub-Finsler Lie groups.
The next result follows from \cite[Proposition 7.1.9]{libroenrico} and the proof therein.

\begin{proposition}
\label{prop:submetries-and-quotients}
    Let $(G,V,\|\cdot\|)$ be a sub-Finsler Lie group and let $H$ be a Lie group. Let $\pi:G\to H$ be a surjective Lie homomorphism.
    On $H$, define the polarization 
    \begin{equation}
        \label{eq:def-polarization-down}
        W :=\df\pi_e(V),
    \end{equation}
and for $w\in W$ set
\begin{eqnarray}
 \|{w}\| : =\inf\{\|v\|:  v\in V, \df\pi_e (v)=w \}.
 \label{eq:def-defnormadown}
\end{eqnarray}
Then $W$ is bracket generating and the map $\pi:G\to H$ between the two sub-Finsler Lie groups is a submetry. 
Moreover, for every interval $I\subseteq \R$, and for every $\gamma:I\to H$ horizontal with control $u$, there exists $\tilde{\gamma}:I\to G$ horizontal with control $ \tilde{u}$ satisfying 
\begin{eqnarray}
    \label{eq:being-lift}
    \pi\circ \tilde{\gamma}&=&\gamma,\\
    \label{preserving-length}
    \|\tilde{u}(t)\|&=&\|u(t)\|, \quad \text{for a.e. } t\in I.
\end{eqnarray}
\end{proposition}

We are now ready to show the main result of this section:
normal curves are lifted to normal curves via sub-Finsler submetries.

\begin{proposition}[Lift of normal curves]\label{prop:lift-normal-curves}
   Let $(G,V,\|\cdot\|)$ be a sub-Finsler Lie group and $\pi:G\to H$ a surjective Lie homomorphism.
 Then there exists a unique sub-Finsler structure on $H$ that makes $\pi$ a submetry;
   explicitly, the sub-Finsler structure on $H$ is the one given by \eqref{eq:def-polarization-down}-\eqref{eq:def-defnormadown}.
   
   Moreover, normal curves lift via $\pi$:
   If $\gamma:I\to H$ is a normal curve with associated covector $\lambda$,
   then every lift $\tilde{\gamma}: I \to G$ (as in \cref{prop:submetries-and-quotients}) is a normal curve with associated covector $\lambda\circ \df \pi_{1}$.
\end{proposition}
\begin{proof}
 The sub-Finsler structure $(W,\|\cdot\|)$ given by \eqref{eq:def-polarization-down}-\eqref{eq:def-defnormadown} makes $\pi$ into a submetry.
 Trivially, there is at most one metric structure on $H$ for which $\pi$ is a submetry.
 
 We now show the second part of the proposition.
    Let $u$ be the control of $\gamma$ and let $\tilde{\gamma}:I\to G$ be a horizontal curve with control $\tilde{u}$ satisfying Equations \eqref{eq:being-lift}-\eqref{preserving-length}. Since the map $\df \pi_e :V\to W $ is a submetry
    , by \cref{cor:lift-subdifferential}, Equation \eqref{preserving-length}, and the normal inclusion \eqref{eq:normal_equation} for $\gamma$, we get 
    \begin{equation}
        \label{eq:pluto1}
        \lambda \Ad_{\gamma(t)}\circ \df \pi_e \in \partial_{\tilde{u}(t)}\E,\quad  \text{for a.e. $t\in I$}.
    \end{equation}
     Moreover, being $\pi:G\to H$ a homomorphism of Lie groups, by \eqref{eq:being-lift} we get 
     \begin{equation}
         \label{eq:pluto2}
         \lambda \Ad_{\gamma}\circ \df \pi_e= \lambda \circ \df \pi_e \circ \Ad_{\tilde{\gamma}}.
     \end{equation}
By Equations \eqref{eq:pluto1} and \eqref{eq:pluto2} it follows that $\tilde{\gamma}$ is normal with covector $\lambda\circ\df \pi_e$.
\end{proof}

\section{Branching of normal curves}
\label{sec:branching}
In this section, we shall prove that on some Finsler Lie groups, there is no possibility of branching of normal curves; see \cref{prop:Finsler-smooth}. Instead, for some (strongly convex) Finsler Lie groups, there is branching; see \cref{them:main-branching-of-normal}. 
We begin by clarifying the notion of branching.

\begin{definition}
    \label{def:branching}
 Let $G$ be a sub-Finsler Lie group. 
 We say that \emph{there is branching of normal curves} in $G$ if there exists two different normal curves $\gamma_1,\gamma_2:I\to G$, where $I\subseteq \R$ is an interval, and a nonempty interval $J\subseteq I$, such that $\gamma_1|_J=\gamma_2|_J$.
\end{definition}

We remark that there are smooth sub-Ri\-em\-man\-ni\-an manifolds that admit branching of normal curves, see \cite{zbMATH07286389};
however, those examples are neither isometrically homogeneous nor analytic.
Sub-Riemannian Lie groups, and more generally analytic sub-Riemannian manifolds,
do not admit branching of normal curves, because in this case, the normal equations and their solutions are analytic. 

It is not known whether sub-Finsler Lie groups with smooth, strongly convex norms may admit branching of normal curves. However, in the next \cref{prop:Finsler-smooth}, we can exclude branching if, in addition, the norm is Finsler.
For non-differentiable strongly convex norms, instead, we will present examples of branching in \cref{them:main-branching-of-normal}.

\subsection{Absence of branching}

Regarding the next result, recall from \cref{cor:uniqueness_sol_strongly_convex} that normal curves in sub-Finsler Lie groups endowed with strongly convex norms are of class $C^{1,1}$, hence it makes sense to consider pointwise derivatives.

\begin{proposition}
\label{prop:Finsler-smooth}
    Let $(G,\g,\|\cdot\|)$ be a Finsler Lie group with a strongly convex norm.
    Let $\gamma_1,\gamma_2:[0,1]\to G$ be normal curves. 
    If there exists $t\in [0,1]$ such that $\gamma_1(t)=\gamma_2(t)$, $\dot\gamma_1(t)=\dot\gamma_2(t)$ and $\|\cdot\|$ is differentiable at $\df L_{\gamma_1(t)}^{-1}\dot \gamma_1(t)$, then $\gamma_1=\gamma_2$.
     In particular, in every Finsler Lie group with strongly convex norm that is differentiable outside 0, there is no branching of normal curves.
\end{proposition}
\begin{proof}
 Denote by $\E:v\mapsto\frac{\|v\|^2}{2}$ the energy of $\|\cdot\|$.
 If $\|\cdot\|$ is differentiable at some $v\neq0$, 
 then $\E$ is differentiable at $v$.
 
    Denote with $u_i$ the control of $\gamma_i$, for $i=1,2$.
    Let $\lambda_1$ and $\lambda_2$ be covectors associated to $\gamma_1$ and $\gamma_2$, respectively. 
    By assumption, there exists $t\in [0,1]$ such that $\gamma_1$ and $\gamma_2$ are differentiable at $t$, $p:=\gamma_1(t)=\gamma_2(t)$, and $u:=u_1(t)=u_2(t)$. 
   Being the metric structure Finsler, its polarization is the Lie algebra $\g$ of $G$. From the normal inclusion~\eqref{eq:normal_equation} it follows that
 \begin{equation*}
     \lambda_1\circ \Ad_{p}|_\g, \  \lambda_2\circ \Ad_{p}|_\g\in \partial_{u}\E.
 \end{equation*}
    Being the energy convex and differentiable at $u$, by \eqref{eq:subdifferential-gateaux-diff} we get
    \begin{equation*}
        \lambda_1\circ \Ad_{p} = \df_u\E = \lambda_2\circ \Ad_{p}.
    \end{equation*}
 Since $\Ad_p$ is invertible, we infer $\lambda_1=\lambda_2$.
 Since $\E$ is strongly convex, \cref{cor:uniqueness_sol_strongly_convex} gives $\gamma_1=\gamma_2$.
\end{proof}

\subsection{Presence of branching}
\label{sec:branching-proof-thm}
In this section, we give a sufficient condition for branching of normal curves. We stress that, as a consequence of  \eqref{eq:subdifferential-gateaux-diff}, in the following theorem, the assumption that $\partial_Y \E $ contains more than one element implies that $\E$ is not differentiable at $Y$. 

\begin{theorem}
\label{them:main-branching-of-normal}
    Let $(G,V,\|\cdot\|)$ be a sub-Finsler Lie group with energy $\E$. 
    Assume that there exist $Y,X\in V$ such that the interior of $\partial_Y \E $ within $S^*(\|Y\|)$ is non-empty and $\Ad_{\exp(t Y)}X$ is unbounded in $t\in\R$. Then there is branching of normal curves in~$G$.  
\end{theorem}

The proof of \cref{them:main-branching-of-normal} is postponed after a lemma.
\begin{lemma}
 \label{rem:parte-interna}
 Let $V$ be a finite-dimensional vector space and $\|\cdot\|$ a norm on $V$ with energy $\E$.
 If $Y\in V\setminus \{0\}$ and $\eta\in \interior_{S^*(\|Y\|)}(\partial_Y\E)$,
 then $\partial_\eta \E^* =\{Y\}$.
\end{lemma}

\begin{proof}
Since $\eta\in \partial_Y\E$, by \cref{remark:duality-subdifferentials} we have $Y\in \partial_\eta \E^*$.
 Suppose by contradiction that there exists $X\in \partial_\eta \E^* \setminus\{Y\}$. By \eqref{eq:def_tan}, the vectors $X$ and $Y$ are on the same sphere and, moreover, cannot be opposite; hence, they are linearly independent. Consequently, there exists $\xi\in V^*$ such that 
 $\langle \xi , Y \rangle=0$ and 
 $\langle \xi, X \rangle = \|X\|$. 
 Set $\eta_\epsilon := \frac{\|Y\|}{\|\eta+\epsilon\xi\|_*}(\eta+\epsilon\xi) \in S^*(\|Y\|)$, for $\epsilon>0$ sufficiently small. 
 Since 
 \[
 \| \eta+\epsilon\xi \|_*
 \stackrel{\eqref{eq:def-dual-norm}}{\geq} 
 \left\langle \eta , \frac{X}{\|X\|} \right\rangle
  + \epsilon \left\langle \xi ,\frac{X}{\|X\|} \right\rangle
 \stackrel{\eqref{eq:def_tan}}{=} \|\eta\|_*+\epsilon
 \stackrel{\eqref{eq:subdiff-of-energy}}{=} \|Y\|+\epsilon ,
 \]
 we have 
 \begin{equation*}
  \langle \eta_\epsilon , Y \rangle
    \stackrel{\eqref{eq:def_tan}}{=} \frac{\|Y\|^3}{\|\eta+\epsilon\xi\|_*}  
     \leq \frac{\|Y\|^3}{\|Y\|+\epsilon}
     < \|Y\|^2
     \stackrel{}{=} \|Y\|\|\eta_\epsilon\|_*.
 \end{equation*}
    In particular, by \eqref{eq:subdiff-of-energy} it follows that $\eta_\epsilon \in S^*(\|Y\|)\setminus \partial_Y\E$ for all $\epsilon>0$ sufficiently small. We reached a contradiction since $\eta\in \interior_{S^*(\|Y\|)}(\partial_Y\E)$.
\end{proof}

\begin{proof}[Proof of \cref{them:main-branching-of-normal}]
 Denote by $\g$ the Lie algebra of $G$.
 We claim that there exists $\lambda\in\g^*$ such that 
 $\lambda|_V\in \interior_{S^*(\|Y\|)}(\partial_Y\E)$ and $$\R \ni t\mapsto \lambda\Ad_{\exp(tY)}|_V \text{ is unbounded}. $$
 Indeed, by assumption there exists $\lambda_1\in \g^*$ such that $\lambda_1|_V$ is in the interior of $\partial_Y \E$ within $S^*(\|Y\|)$.
 If $t\mapsto \lambda_1(\Ad_{\exp(tY)}X)$ is unbounded in $t\in\R$, set $\lambda:=\lambda_1$. Otherwise, choose $\lambda_2\in \g^*$ with the property that $\lambda_2(\Ad_{\exp(tY)}X)$ is unbounded in $t\in\R$, which exists by assumption, and set $\lambda:=\frac{\|Y\|}{\|(\lambda_1+\alpha\lambda_2)|_V\|_*} (\lambda_1+\alpha\lambda_2)$, with $\alpha>0$ small enough in order to have $\lambda|_V\in \interior_{S^*(\|Y\|)}(\partial_Y\E)$.
   We thus have the covector we looked for.
   
   Let $\gamma:\R \to G$ be a normal curve with covector $\lambda$ such that $\gamma(0)$ is the unit element of $G$,
   which exists by \cref{prop:solution-normal-exists}.
   By \cref{prop:normal-PBAL}, the control $u$ of $\gamma$ has norm constantly equal to $\|\lambda|_V\|_*=\|Y\|$. By \eqref{eq:normal_equation} and \cref{prop:subdifferentials-of-energy} we have that $\lambda\Ad_{\gamma(t)}|_V\in \partial_{u(t)}\E\subseteq S^*(\|u(t)\|)=S^*(\|Y\|)$, for almost every $t\in \R$.

Set 
\begin{equation*}
    \bar{t}:=\sup \{t \in [0,+\infty) : \forall\tau \in [0,t], \gamma(\tau)=\exp(\tau Y) \}.
\end{equation*}
Notice that, on the one hand, we just proved that $t\mapsto \lambda\Ad_{\gamma(t)}|_V$ is bounded. On the other hand, we constructed $\lambda$ so that $t\mapsto\lambda\Ad_{\exp(tY)}$ is unbounded. We deduce that $\bar t<\infty$.

We claim $\bar t>0$.
Since $\lambda|_V\in \interior_{S^*(\|Y\|)}(\partial_Y\E)$ and $\gamma$ is continuous, there exists $\epsilon>0$ such that $\lambda\Ad_{\gamma(t)}|_V\in \interior_{S^*(\|Y\|)}(\partial_Y\E)$ for all $t\in [0,\epsilon]$.
By \cref{rem:parte-interna}, there holds
$\partial_{\lambda\Ad_{\gamma(t)}|_V}\E^* = \{Y\}$ for all $t\in [0,\epsilon]$.
Since $\gamma$ is a normal curve, it satisfies~\eqref{eq:normal-differential-inclusion}, hence
 we get $\gamma(t)=\exp(tY)$ for all $t\in[0,\epsilon]$. 
 Thus $\bar{t}\ge\epsilon>0$, and the claim is proven. 

  
For $\beta\in [0,\bar{t})$, let $\gamma_\beta:\R\to G$ be the curve $\gamma_\beta(t):=  \gamma(\beta)^{-1}\gamma(t+\beta)$. The curve $\gamma_\beta$ is normal with associated covector $\lambda_\beta:=\lambda\Ad_{\exp(\beta Y)}$.
We have
   \begin{eqnarray*}
     \gamma_\beta(t)=\exp(tY) = \gamma(t),\quad \forall  t\in [0,\bar{t}-\beta].
   \end{eqnarray*}
By the definition of $\bar t$, there is $\tau\in(\bar t-\beta,\bar t)$ with $\gamma(\tau+\beta) \neq \exp((\tau+\beta)Y)$ and thus $\gamma_\beta(\tau) \neq \exp(\tau Y) = \gamma(\tau)$. 
%
%
Therefore, the curves $\gamma_\beta$, for $\beta\in [0,\bar{t})$, are normal curves that all coincide with the one-parameter subgroup $t\mapsto \exp(tY)$ on $[0,\bar{t}-\beta]$, and are all different curves.
\end{proof}


\cref{thm:if-condition-then-exists-branching} follows from \cref{them:main-branching-of-normal} as soon as we can show the existence of a norm with the property that $\interior_{S^*(\|Y\|)}{\partial_Y\E}\neq \emptyset$.
We do it in \cref{prop:esistenza-norma-strong-convex}. 


\begin{lemma}
 \label{prop:esistenza-norma-strong-convex}
 Let $V$ be a finite-dimensional vector space and $v\in V\setminus\{0\}$. 
 Then there exists a strongly convex norm $\|\cdot\|$ on $V$ 
 such that $\interior_{S^*({\|v\|})}{\partial_v\E}\neq \emptyset$,
 where $\E$ is the energy of $\|\cdot\|$.
\end{lemma}
\begin{proof}
    Without loss of generality, we assume $V=\R^n$ and $v$ to be the first vector $e_1$ of the canonical base $e_1,...,e_n$ of $\R^n$. We consider the norm
    \begin{equation*}
        \|\cdot\|:=\sqrt{\|\cdot\|_{L^1}^2+\|\cdot\|_{L^2}^2 }.
    \end{equation*}
On the one hand, the energy $\E_{\|\cdot\|}$ is strongly convex being the sum of the convex function $\E_{\|\cdot\|_{L^1}}$ and the strongly convex function $\E_{\|\cdot\|_{L^2}}$. On the other hand, by \eqref{eq:sum-of-subdifferentials} we have
\begin{equation*}
    \partial_{e_1}E_{\|\cdot\|}= \partial_{e_1}\E_{\|\cdot\|_{L^1}} + \partial_{e_1}\E_{\|\cdot\|_{L^2}},
\end{equation*}
    hence $\partial_{e_1}E_{\|\cdot\|}$ contains an open subset of an hyperplane since $\partial_{e_1}\E_{\|\cdot\|_{L^1}}$ does. This implies that $\interior_{S^*({\|v\|})}{\partial_v\E}$ is non-empty.
\end{proof}

\subsection{Example: branching in a $2D$-group with strongly convex norm}
\label{sec:example-2-dim-non-ab-group}
In this section, we present examples of branching. The condition from \cref{them:main-branching-of-normal} that there exists $X,Y\in \g$ such that the set $\{ \Ad_{\exp(t Y)}X : t\in\R\}$ is unbounded is satisfied by several Lie algebras, such as, for example, all nilpotent non-abelian Lie algebras. We stress that for the latter Lie algebras the function $t\mapsto \Ad_{\exp(t Y)}X$ is polynomial in $t$.

\begin{figure}[H]
    \centering
    \begin{tikzpicture}[scale=2]
        \draw[->] (-1.2,0) -- (1.2,0) node[below] {\footnotesize $x$};
        \draw[->] (0,-1.3) -- (0,1.3) node[left] {\footnotesize $y$};

        \draw[thick, domain=-1:1, samples=100, smooth, variable=\y]
            plot ({(\y*\y)/2 - 0.5}, {\y});

        \draw[thick, domain=-1:1, samples=100, smooth, variable=\y]
            plot ({-((\y*\y)/2 - 0.5)}, {\y});

        \filldraw (0, 1) circle (0.02) node[right, yshift=0.2cm] {\footnotesize $(0,1)$};
        \filldraw (0, -1) circle (0.02) node[right] {\footnotesize $(0,-1)$};
        \filldraw (0.5, 0) circle (0.02) node[above right] {\footnotesize $(\tfrac{1}{2},0)$};
        \filldraw (-0.5, 0) circle (0.02) node[above left] {\footnotesize $(-\tfrac{1}{2},0)$};

        \draw[blue, thick] (-0.6, 1.3) -- (0.6, 0.7) node[right] {\footnotesize $\lambda=(\frac{1}{2},1)$};

        \draw[red, thick] (-0.6, 1.2) -- (0.6, 0.8) node[right, yshift=0.2cm] {\footnotesize $\lambda'=(\frac{1}{3},1)$};
    \end{tikzpicture}
    \caption{The unit sphere of the norm $\|\cdot\|$ is the non-smooth strongly convex set $\{(x,y)\in\mathbb{R}^2 : |x| = \frac{y^2}{2} - \frac{1}{2}\}$. The blue and red lines correspond to two different covectors $\lambda$ and $\lambda'$ defining supporting hyperplanes of the sphere.}
    \label{figure:covectors}
\end{figure}
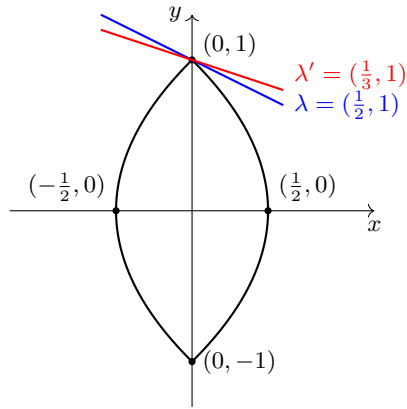

\begin{figure}[H]
    \centering
    \includegraphics[width=0.5\linewidth]{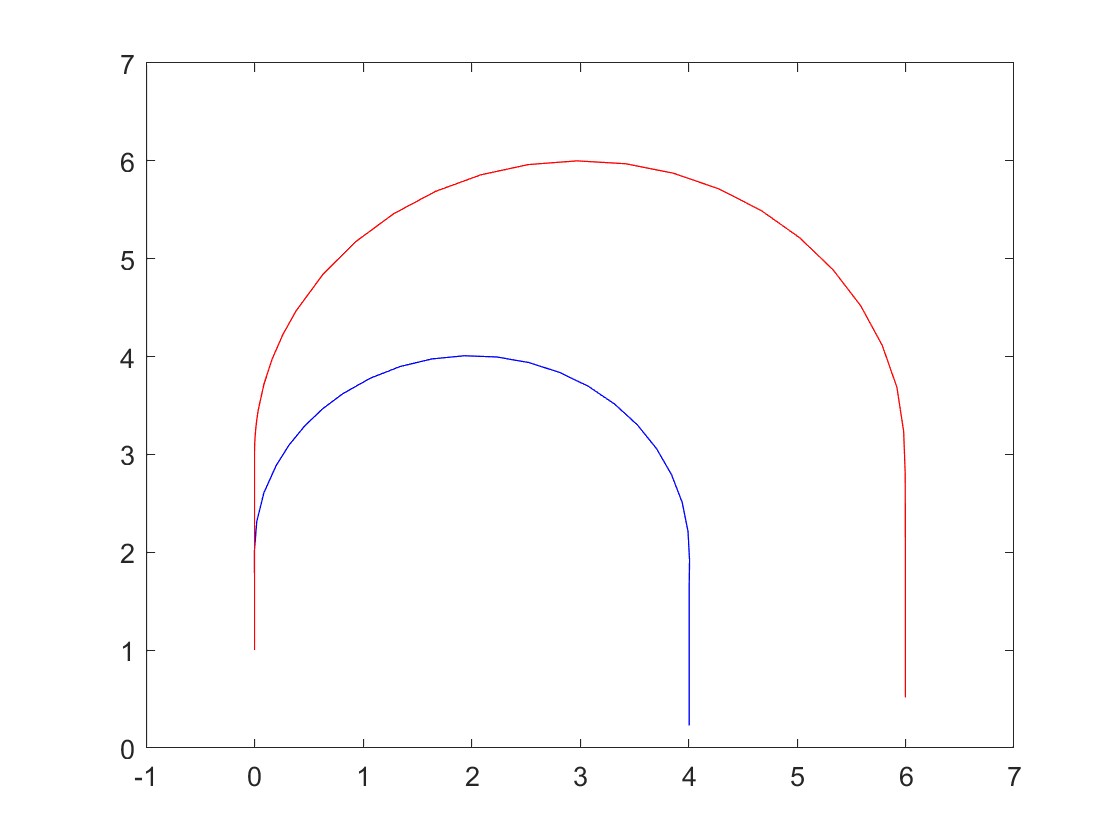}
    \caption{The numerical integration of the normal equation for the covectors $\lambda$ and $\lambda'$.
    In blue, the normal curve with associated covector $\lambda$ starting at $(0,1)$. In red, the normal curve with associated covector $\lambda'$ starting at $(0,1)$.}
    \label{fig:num_int}
\end{figure}

We shall provide another concrete example of branching in the $2$-dim\-en\-sion\-al non-abelian Lie algebra. Namely, let $G$ be the group $\R\rtimes \R_+$ with product law 
\begin{equation*}
    (x,t)\cdot(y,s):=(x+ty, ts),\quad \forall x,y\in\R,s,t\in \R_+.
\end{equation*}
In the Lie algebra $\R^2$ of $G$ consider the norm
\begin{equation*}
    \|(x,y)\|:= |x|+\sqrt{x^2+y^2},\quad \forall (x,y)\in\R^2.
\end{equation*}
Let us consider $\lambda:=(\frac{1}{2},1)$ and $\lambda':=(\frac{1}{3},1)$, as covectors at the identity point $(0,1)$, see \cref{figure:covectors}. The normal equations for the curves with covectors $\lambda$ and $\lambda'$ are 
\begin{equation*}
\begin{cases}
    \dot \gamma= 
\gamma_2F\left(\frac{\gamma_2}{2}, -\frac{\gamma_1}{2}+1\right),\\
    \gamma(0)=(0,1),
\end{cases}
\end{equation*}
and
\begin{equation*}
\begin{cases}
    \dot \gamma= 
\gamma_2F\left(\frac{\gamma_2}{3}, -\frac{\gamma_1}{3}+1\right),\\
    \gamma(0)=(0,1),
\end{cases}
\end{equation*}
where $F:\R^2\to \R^2$ is defined by setting, for all  $\eta\in\R^2$,
\begin{equation*}
    F(\eta_1,\eta_2):=\begin{cases}
        (0,\mathrm{sign}(\eta_2)),\quad \text{ if } |\frac{\eta_2}{\eta_1}|\geq 1;\\
        \\
        \mathrm{sign}(\eta_1)\left(\frac{1}{2}-\frac{1}{2}\left(\frac{\eta_2}{\eta_1}\right)^2, \frac{\eta_2}{\eta_1}\right),\quad \text{ if } |\frac{\eta_2}{\eta_1}|\leq 1.
    \end{cases}
\end{equation*}
We refer to \cref{fig:num_int} for a numerical integration of the associated normal curves. We stress that the two normal curves coincide with the same one-parameter subgroup on the intervals $[0,\log(2)]$ and $[0,\log(3)]$ respectively. Afterwards, they branch apart. 

\subsection{Example: the special orthogonal group}\label{subs68a6d100}
Clearly, in abelian groups with strictly convex norms, there is no branching of normal curves. In this section, we present a non-abelian example where the energy is strongly convex and not smooth, and still, there is no branching. 
As an example, we consider a Finsler norm on the special orthogonal group $SO(3)$ with the property that there exists $Y\in \mathfrak{so}(3)$ such that the interior of $\partial_Y \E $ within $S^*(\|Y\|)$ is non-empty. However, in this particular (sub-)Finsler group, there is no branching of normal geodesics. This shows that some hypotheses on the adjoint map, e.g., the unboundedness of the adjoint in \cref{them:main-branching-of-normal}, are necessary if we want to ensure the presence of branching. We consider the base of $\mathfrak{so}(3)$ given by
\begin{eqnarray*}
    e_1 := 
\begin{pmatrix}
0 & 1 & 0 \\
-1 & 0 & 0 \\
0 & 0 & 0
\end{pmatrix}, \quad
e_2 := 
\begin{pmatrix}
0 & 0 & 1 \\
0 & 0 & 0 \\
-1 & 0 & 0
\end{pmatrix}, \quad
e_3 := 
\begin{pmatrix}
0 & 0 & 0 \\
0 & 0 & 1 \\
0 & -1 & 0
\end{pmatrix}.
\end{eqnarray*}

On $\mathfrak{so}(3)$ we define the norm, for $(x_1,x_2,x_3)\in \R^3$,
\begin{eqnarray*}
    \left\| x_1 e_1 + x_2 e_2 + x_3 e_3 \right\| := \sqrt{x_2^2 + x_3^2} + \sqrt{x_1^2 + x_2^2 + x_3^2}.
\end{eqnarray*}
We check that $\partial_{e_1}\E$ has non-empty interior within the two dimensional sphere $S^*$.
Indeed, we consider the dual base $e_1^*,e_2^*,e_3^*$ of $e_1,e_2,e_3$ and for all $\alpha,\beta\in (-\frac12,\frac{1}{2})$ we define $\lambda_{\alpha,\beta}:=e_1^*+\alpha e_2^*+\beta e_3^*$.
Then, we have
$\lambda_{\alpha,\beta}(e_1)=1$ and, for all $x:=x_1e_1+x_2e_2+x_3e_3$, with $(x_1,x_2,x_3)\in \R^3$, there holds 
\begin{align*}
    \lambda_{\alpha,\beta}(x)&\leq |x_1|+|\alpha| |x_2|+|\beta| |x_3|\\&\leq |x_1|+\frac12 |x_2|+\frac12 |x_3|\\
    &\leq |x_1|+\sqrt{x_2^2+x_3^2} \leq \|x\|.
\end{align*}
Hence, by \eqref{eq:subdiff-of-energy}, we have $\{\lambda_{\alpha,\beta} \mid \alpha,\beta\in (-\frac12,\frac{1}{2})\}\subseteq \partial_{e_1}\E$.

The norm $\|\cdot\|$ is strongly convex and smooth outside the set $\R e_1$. Hence, by \cref{prop:Finsler-smooth}, we have that if two normal curves have control that takes value outside of $\R e_1$ and coincide on an open interval, then they coincide everywhere. To show that there is no branching of normal curves, we only need to show that there is no branching between normal curves whose control takes value in $\R e_1$ at some time. To do so, we will show that if a normal curve has control in $\R e_1$ at some time, then it is a re-parametrization of the one-parameter subgroup $t\mapsto \exp(t e_1)$. 
A simple computation shows that 
\begin{align*}
\Ad_{\exp(te_1)} e_1&=e_1&,\\
\Ad_{\exp(te_1)} e_2 &= 
 \begin{pmatrix}
 \cos t & \sin t & 0 \\
 -\sin t & \cos t & 0 \\
 0 & 0 & 1
 \end{pmatrix}
 \begin{pmatrix}
 0 & 0 & 1 \\
 0 & 0 & 0 \\
 -1 & 0 & 0
 \end{pmatrix}
 \begin{pmatrix}
 \cos t & -\sin t & 0 \\
 \sin t & \cos t & 0 \\
 0 & 0 & 1
 \end{pmatrix} 
 \\
 &= \cos t \, e_2 + \sin t \, e_3,
\end{align*}

\begin{align*}
\Ad_{\exp(te_1)}  e_3 &= 
\begin{pmatrix}
\cos t & \sin t & 0 \\
-\sin t & \cos t & 0 \\
0 & 0 & 1
\end{pmatrix}
\begin{pmatrix}
0 & 0 & 0 \\
0 & 0 & 1 \\
0 & -1 & 0
\end{pmatrix}
\begin{pmatrix}
\cos t & -\sin t & 0 \\
\sin t & \cos t & 0 \\
0 & 0 & 1
\end{pmatrix} \\
&= \sin t \, e_2 + \cos t \, e_3.
\end{align*}
In particular,

\begin{equation}
    \label{eq:norm-constant-along-Ad-ops}
    \| \Ad_{\exp(t e_1)}v\|=\|v\|, \quad \forall v\in \mathfrak{so}(3), \forall t\in\R.
\end{equation}
Fix $\alpha\in \R\setminus \{0\}$. We claim that for all $\lambda\in \mathfrak{so}(3)^*$ we have 
\begin{equation}
    \label{claim:covecto-ops-so3}
    \lambda\in \partial_{\alpha e_1}\E \implies \lambda \Ad_{\exp(t \alpha e_1)}\in \partial_{\alpha e_1}\E, \quad \forall t\in\R.
\end{equation}
Indeed, if $\lambda\in \partial_{\alpha e_1}\E$, then
\begin{equation*}
    \lambda\Ad_{\exp(t \alpha e_1)}(\alpha e_1)=\lambda(\alpha e_1)\stackrel{\eqref{eq:subdiff-of-energy} }{=} \alpha^2,
\end{equation*}
and moreover, by \eqref{eq:norm-constant-along-Ad-ops} the map $\Ad_{\exp(t \alpha e_1)}$ is an isometry,
hence $$\|\lambda \Ad_{\exp(t \alpha e_1)}\|_*=\|\lambda\|_*=\|\alpha e_1\|=\alpha.$$ Consequently, the claim \eqref{claim:covecto-ops-so3} follows by \eqref{eq:subdiff-of-energy}.

As a consequence of \eqref{claim:covecto-ops-so3} and the normal inclusion \eqref{eq:normal_equation}, we obtained 
\begin{equation}
\label{eq:exp_e_1-ops-for-all-covectors}
     \lambda\in \partial_{\alpha e_1}\E \implies t\in\R \mapsto \exp(t\alpha e_1) \text{ is normal with covector } \lambda. 
\end{equation}

Let $\gamma:[0,1]\to SO(3)$ be a normal curve for the above Finsler structure with associated covector $\lambda\in \mathfrak{so}(3)$.
Assume that $\gamma$ has control equal to $\alpha e_1$ at some time $t_0\in \R$, for some $\alpha\in\R \setminus \{0\}$. Up to translations we assume $t_0=0$ and $\gamma(0)=1_G$. Then, by the normal inclusion \eqref{eq:normal_equation}, we have $\lambda\in \partial_{\alpha e_1}\E$. By \eqref{eq:exp_e_1-ops-for-all-covectors} and by uniqueness of the solution of the normal equation for strongly convex norms, see \cref{cor:uniqueness_sol_strongly_convex}, we get $\gamma(t)=\exp(t\alpha e_1)$ for all $t\in[0,1]$. Therefore, in this Finsler Lie group, there is no branching of normal curves, even if the norm is not smooth and the group is not abelian.

\section{Sub-Finsler Lie groups with polyhedral norms}
\label{sec:polytope}
\subsection{Polyhedral norms}

\begin{definition}\label{def:polytope-norm}
Let $V$ be a finite-dimensional vector space. 
We say that $P\subseteq V$ is a \textit{polyhedron} if it is the convex hull of finitely many points. A norm on $V$ is a \emph{polyhedral norm} if the closed unit ball is a polyhedron.
\end{definition}
For every polyhedron $P$ that satisfies $P=-P$ and $0\in \interior(P)$, there exists a unique norm $\|\cdot\|_P$ whose unit ball coincides with $P$, and there holds
\begin{equation*}
    \|v\|_P:=\min \{\alpha\in \R : \frac{1}{\alpha}v\in P\}.
\end{equation*}
We call $\|\cdot\|_P$ the \emph{norm associated to} $P$.
\begin{remark}
\label{rem:poly-norms}
    Let $V$ be a finite-dimensional vector space and $P\subseteq V$ a polyhedron. Then there exists a minimal (with respect to the inclusion) finite subset $\Lambda_P\subseteq V^*$ with the property that
    \begin{equation*}
        P = \{v\in V : \lambda(v)\leq 1, \ \lambda\in \Lambda_P\}.
    \end{equation*}
    If $P$ is the unit ball of a polyedral norm $\|\cdot\|_P$, then 
    \begin{equation*}
        \|v\|_P=\max\{\lambda(v) \ : \ \lambda\in\Lambda_P\}.
    \end{equation*}
    In this case, the convex hull of $\Lambda_P$ in $V^*$ is the unit ball of the dual norm $\|\cdot\|_*$.
    In particular, the dual norm of a polyhedral norm is polyhedral,
    and both the unit sphere $S$ and the dual sphere $S^*$ are simplicial complexes;
    see also \cite[Chapter 19]{Clarke-book}.
\end{remark}

\begin{definition}[Face]
 Let $K\subseteq V$ be a compact convex set. A set $F\subset K$ is a \emph{face} of $\partial K$ if there exists $\eta\in V^*$ such that
 \begin{equation*}
     F=F_\eta:=\{v\in K \mid \langle\eta,v\rangle= \max_{w\in K}\langle\eta,w\rangle.
 \end{equation*}
\end{definition}

\begin{remark}
    By \cref{prop:subdifferentials-of-energy}, for all $\eta\in V^*$ with $\|\eta\|_*=1$, we have
\begin{equation}
    \label{eq:faccia-associata-eta}
    \partial_\eta\E^*=\{v\in V\mid \|v\|=\langle\eta,v\rangle=1\}\subseteq S.
\end{equation}
Hence, the set $\partial_\eta\E^*\subseteq S$ is the face of $S$ determined by $\eta$.
    
\end{remark}

When studying polyhedral norms, we will be interested in the following set of covectors:

\begin{definition}
\label{def:tan-e-star}
   Let $V$ be a vector space and $E:V\to\R$ the energy of a norm on $V$.
   For $\eta\in V^*$, the \emph{\nomedascegliere} of $\eta$ is the set 
    \begin{equation*}
            \starr(\eta):=\{\xi\in V^* : \partial_\xi \E^* \subseteq \partial_\eta \E^*\}.
    \end{equation*}
\end{definition}

\begin{remark}
Let $\eta\in S^*$. By \cref{prop:subdifferentials-of-energy}, we have that $\xi\in \starr(\eta)$ if and only if $\|\xi\|_*=1$ and for all $v\in S$ with $\langle \xi,v\rangle=1$ we have $\langle \eta,v\rangle=1$. Hence, again by \cref{prop:subdifferentials-of-energy}, we have that $\xi\in \starr(\eta)$ if and only if, for all $v\in S$ with $\xi\in\partial_v\E$ we have $\eta\in\partial_v\E$. Consequently, 
\begin{equation}
\label{eq:star-as-union}
    \starr(\eta)= S^* \setminus \bigcup_{\substack{v \in S \\[2pt] \eta \notin \partial_v \E}} \partial_v \E,
\end{equation}
that is, the set $\starr(\eta)$ is the complementary in $S^*$ of all the faces that do not contain $\eta$.
In particular, in the case the norm is polyhedral, if $\eta$ is a vertex of the simplicial complex $S^*$, then $\starr(\eta)$ coincides with the open star of $\eta$ within the simplicial structure. For the definition of the latter, see \cite[Section 0.2]{stillwell2012classical}.
\end{remark}


Polyhedral norms can be characterized in terms of the sets $\partial \E^* $ and $\starr$ defined in \cref{def:tan-e-star}. 
\begin{proposition}
\label{prop:polytope-norm-equivalences}
    Let $V$ be a finite-dimensional vector space and $\|\cdot\|$ a norm on $V$. The following are equivalent.
    \begin{enumerate}[label=(\roman*)]
        \item
        The norm $\|\cdot\|$ is a polyhedral norm;
        \item The set $\{\partial_\eta \E^*  : \eta\in S^*\}$ is finite;
        \item There exist $\eta_1,...,\eta_n\in S^*$ such that $\{\starr(\eta_i)\}_{i\in \{1,...,n\}}$ is an open covering of $S^*$. 
    \end{enumerate}
\end{proposition}
\begin{proof}
By \cite[Theorem 19.1]{rockafellar2015convex}, a convex set is a polyhedron if and only if it has a finite number of faces; thus, properties (i) and (ii) are equivalent. 
  
  We next prove that (iii) implies (i). 
  By \eqref{eq:faccia-associata-eta}, we have $S=\cup_{\eta\in S^*} \partial_\eta\E^*$. By (iii), for all $\eta\in S^*$ there exists $i\in\{1,...,n\}$ such that $\eta\in \starr(\eta_i)$, that is, we have $\partial_\eta\E^* \subseteq \partial_{\eta_i}\E^*$. Hence, we get $S=\cup_{i=1}^n \partial_{\eta_i}\E^*$. We showed that the sphere is a union of a finite number of faces, thus the norm is polyhedral.
  
  Finally, we prove that (ii) implies (iii). 
  Assume that there exist $\eta_1,...,\eta_n$ such that
  \begin{equation}
  \label{eq:claim-in-remark-poliedri-star}
      \{\partial_{\eta} \E^*  : \eta\in S^*\}=
      \{\partial_{\eta_i} \E^*  : i\in \{1,...,n\}\}.
  \end{equation}
Then by definition of star we deduce that $\{\starr(\eta_i)\}_{i\in\{1,...,n\}}$ is a covering of $S^*$. We claim that $\starr(\eta_i)$ is open for all $i \in \{1, \ldots, n\}$. 
Indeed, since the norm $\|\cdot\|$ is polyhedral, its dual norm $\|\cdot\|_*$ is also polyhedral, see \cite[Chapter~19]{rockafellar2015convex}. Hence, the set $\{\partial_v \E : v \in S\}$ is finite. 
Consequently, by~\eqref{eq:star-as-union}, we deduce that for each $i \in \{1, \ldots, n\}$, the set $S \setminus \starr(\eta_i)$ is a finite union of closed sets.
  
\end{proof}

\subsection{Normal curves in sub-Finsler groups with polyhedral norms.}\label{subs689c584c}
This section is devoted to the proof of \cref{thm:main-polytope-norms-projection2} and its consequences. We start by proving that on sub-Finsler Lie groups with polyhedral norm, the control of normal curves takes values in a face on intervals of small size.

\begin{proof}[Proof of \cref{thm:main-polytope-norms-projection2}] 
 Let $(G,V,\|\cdot\|)$ be a sub-Finsler Lie group where $\|\cdot\|$ is a polyhedral norm.
 We denote by $\E:V\to\R$ the energy function  and by $\g$ the Lie algebra of $G$.
 Fix an auxiliary norm $N:\g\to\R$.
Since all norms on finite-dimensional vector spaces are bi-Lipschitz equivalent and since \cref{thm:main-polytope-norms-projection2}.(ii) is invariant under bi-Lipschitz change of the norm $N$, 
we can assume $N(v)\leq \|v\|$ for all $v\in V$. Similarly, we take $N^*$ to be the dual norm of $N$.

    By \cref{prop:polytope-norm-equivalences}, there exist $\eta_1,...,\eta_n\in S^*$ such that $\{\starr(\eta_i)\}_{i=1}^n$ is an open covering of $S^*$. 
    Let $\delta>0$ be the Lebesgue number of the covering $\{\starr(\eta_i)\}_{i=1}^n$ with respect to $\|\cdot\|_*$.
    Hence, for all curves $\xi:[0,1]\to S^*$, if $\mathrm{Length}_{\|\cdot\|_*}(\xi)\leq \delta$, then there exists $i\in\{1,...,n\}$ such that $\xi(t)\in \starr(\eta_i)$ for all $t\in [0,1]$. Observe that, for all $r>0$, the family of sets $\{\starr(r\eta_i)\}_{i=1}^n$ is an open covering of $S^*(r)$ with Lebesgue number~$r\delta$.

    Let $\gamma:[0,1]\to G$ be a normal curve such that $\gamma(0)=1_G$ with control $u$ and associated covector $\lambda$. By \cref{prop:normal-PBAL}, the curve $\gamma$ is parametrized by multiple of arclength. Denote with $r$ the speed of $\gamma$.
    
    Fix $t_1,t_2\in [0,1]$, with $t_1<t_2<t_1+\frac{\delta}{N^*(\lambda)M\left(r\right)}$, where $M$ is defined by \eqref{eq:def-M} and $N^*$ is the dual norm of $N$.
     By the normal inclusion \eqref{eq:normal_equation} and \cref{prop:subdifferentials-of-energy} we have  $\lambda\Ad_{\gamma(t)}|_V\in S^*(r)$ for all $t\in [0,1]$. Moreover, for almost every $t\in [0,1]$, we use that $\Ad_*=\ad$ and bound
    \begin{eqnarray*}    
    \left\|\left(\lambda\Ad_{\gamma(t)}|_V\right)'\right\|_*
    &\stackrel{}{=}&\|\lambda\Ad_{\gamma(t)} \ad_{u(t)}|_V\|_*\\ 
 &\stackrel{\eqref{eq:def-dual-norm}}{=}& \max\{\lambda\Ad_{\gamma(t)} \ad_{u(t)}(Y) : Y\in V,\|Y\|=1 \}\\
 &\leq& \|u(t)\|\max \left\{ \lambda\Ad_{\gamma(t)}\ad_X(Y) :  \begin{array}{c}
     {X,Y\in V,}\\ {N(X)=N(Y)=1}
 \end{array} \right\}\\
 &\stackrel{\eqref{eq:def-dual-norm}}{\leq}& r N^*(\lambda) \max\left\{ N(\Ad_{\gamma(t)}\ad_X(Y)) : \begin{array}{c}
     {X,Y\in V,}\\{N(X)=N(Y)=1}
 \end{array} \right\}\\
   &\stackrel{\eqref{eq:def-M}}{\leq}& rN^*(\lambda)M(r),
   \end{eqnarray*}
where in the last inequality we used that $\gamma(t)\in B(1_G,r)$.
Consequently,
\begin{eqnarray*}
    \mathrm{Length}_{\|\cdot\|_*}((\lambda\Ad_{\gamma}|_V)|_{[t_1,t_2]})&\leq& |t_1-t_2|rN^*(\lambda)M(r)\\
    &\leq& \frac{\delta rN^*(\lambda)M(r)}{N^*(\lambda)M\left(r\right)}= \delta r .
\end{eqnarray*}

Since the covering $\{\starr(r\eta_i)\}_{i=1}^n$ of $S^*(r)$ has Lebesgue number $r\delta$, there exists $i\in\{1,...,n\}$ such that
\begin{equation}
\label{eq:Ad-gamma-in-star-eta}
    \lambda\Ad_{\gamma(t)}|_V\in \starr(r\eta_i)
\end{equation}
for all $t\in [t_1,t_2]$.
By \eqref{eq:normal-differential-inclusion}, for all $t\in[t_1,t_2]$ we have
\begin{eqnarray*}
    u(t)\stackrel{\eqref{eq:normal_equation}}{\in} \partial_{\lambda\Ad_{\gamma(t)}|_V} \E^* \subseteq \partial_{r\eta_i} \E^*,
\end{eqnarray*}
where the last inclusion follows by \eqref{eq:Ad-gamma-in-star-eta} and \cref{def:tan-e-star}.
This shows that the control $u|_{[t_1,t_2]}$ takes values only on the face $\partial_{r\eta_i} \E^* $ and proves the theorem.
\end{proof}

We next pass to the case of Finsler groups, where we prove that the controls of normal curves shorter than a certain length take values in a face.
\begin{proof}[Proof of \cref{cor:finsler-compact}]
Let $(G,\mathfrak{g}, \|\cdot\|)$ be a Finsler Lie group, where $\g$ is the Lie algebra of $G$. Consider $N:= \|\cdot\|$ as a norm on $\g$ and $N^*$ its dual norm. Let $\delta$ be the constant coming from \cref{thm:main-polytope-norms-projection2}. Since $\lim_{L\to 0} \frac{\delta}{LM(L)}=+\infty$, there exists $\bar{L}>0$ such that $\frac{\delta}{LM(L)}>1$ for all $L\in (0,\bar L)$.
We claim that, for all $L\in (0,\bar L)$, for every normal curve $\gamma$ of length $L$, the control of $\gamma$ takes value in a face.
Let $\gamma:[a,b]\to G$ be a normal curve of length $L\in (0,\bar L)$ and covector $\lambda\in \g^*$. Up to perform a left-translation and an affine reparametrization, we assume
$\gamma:[0,1]\to G$ and that $\gamma(0)$ is the unit element of $G$. By \cref{prop:normal-PBAL}, we have that $\gamma$ is parametrized with constant speed $L$ and, being $V=\g$, there holds $N^*(\lambda)=L$.
Then, being $\frac{\delta}{LM(L)}>1$, by \cref{thm:main-polytope-norms-projection2} the control of $\gamma$ takes value in a face on $[0,1]$.
\end{proof}
\subsection{Carnot groups with polyhedral norms}
\label{sec:carnot-polyhedral}



In this section, we consider some sub-Finsler Carnot group $G$; we refer to \cite{libroenrico} for a comprehensive introduction to Carnot groups. More generally, we could consider a nilpotent simply connected Lie group $G$, with Lie algebra $\g$ that is polarized by some subspace $V\subseteq \g$ so that $\g=V\oplus [\g,\g]$. On such a $V$, we fix a norm. In this situation, the abelianization $G/[G,G]$ of $G$ admits a unique normed vector space that makes the canonical projection $G\to G/[G,G]$ a submetry (see \cref{def:submetry}). 
Every absolutely continuous curve in $G/[G,G]$ lifts to a horizontal curve in $G$, and the lift is unique up to the choice of an initial point.
Energy minimizing curves in $G/[G,G]$ lift as energy minimizing in $G$.

 Therefore, horizontal curves in $G$ that project to energy minimizers in $G/[G,G]$ are energy minimizers in $G$. The interested reader can find more details in \cite[Section~10.1.1]{libroenrico}.

In sub-Finsler Carnot groups with polyhedral norm, we show that if the control of a curve stays on a face of the polyhedral sphere, then the projection of the curve to the abelianization of the group is a geodesic (up to reparametrization). 

\begin{proof}[Proof of \cref{cor:minimality-normal-curves}]

By \cref{prop:polytope-norm-equivalences}, the set $\{\partial_\eta \E^*  : \eta\in S^*\}$ of the faces of $S$ is finite, that is, we have $\{\partial_\eta \E^*  : \eta\in S^*\}=\{\partial_{\eta_i} \E^*  : i\in 1,...,n\}$, for some $n\in \N$, $\eta_1,...,\eta_n\in S^*$.  
Let $\delta$ be the constant coming from \cref{thm:main-polytope-norms-projection2} and let $\gamma:[0,1]\to G$ be a normal curve with control $u$ and speed $r$ such that $\gamma(0)=1_G$.
Fix $t_1,t_2\in [0,1]$, with $t_1<t_2<t_1+\frac{\delta}{N^*(\lambda)M\left(r\right)}$. 
By \cref{thm:main-polytope-norms-projection2} there exists $i\in \{1,...,n\}$ such that 
\begin{eqnarray}
\label{eq:constant-subdifferential2}
    u(t)\in \partial_{r \eta_i} \E^*, \quad \text{for a.e. } t\in[t_1,t_2].
\end{eqnarray}
The curve $\pi\circ \gamma:[0,1]\to G/[G,G]$ has $u$ as control, which by \eqref{eq:constant-subdifferential2} takes value in a face on the interval $[t_1,t_2]$. Consequently, since $G/[G,G]$ is a normed vector space, the curve $\pi\circ\gamma|_{[t_1,t_2]}$ is length-minimizing by a classical argument, see for example \cite[Lemma 1.7]{creutz2021rigidity}. 
By the discussion before the proof, the curve $\gamma|_{[t_1,t_2]}$ is length-minimizing. Hence, we proved the corollary for $\alpha:=\frac{\delta}{M(1)}$.
\end{proof}

\subsection{A normal extremal that is not locally optimal}
\label{sec:non-locally-optimal-extremal-engel}

In this section, we present a normal curve in a sub-Finsler Lie group that is not locally optimal. 
In fact, we prove a stronger statement: we show the existence of a covector $\lambda$ such that all normal curves with covector $\lambda$ fail to be locally optimal.

Recall that normal curves in sub-Finsler Lie groups with strongly convex norms are locally optimal; see \cref{rem68ad82a9}. Our example will not be strongly convex.

Let  $G$ be the Heisenberg group, that is, the simply connected nilpotent group with Lie algebra $\mathfrak{g}$ spanned by three vectors $X_1,X_2,X_3$ satisfying, as the only nontrivial bracket relation,
\begin{equation*}
    [X_1,X_2]=X_3.
\end{equation*}
On the Lie algebra of $G$, we consider the norm
\begin{equation*}
    \left\|\sum_{i=1}^3 x_iX_i\right\|:=\max\{|x_i| : i\in\{1,2,3\}\},\quad \forall x\in\R^3.
\end{equation*}
The sub-Finsler Lie group $(G,\mathfrak{g},\|\cdot\|)$ has a polyhedral norm,
but it is not Carnot because the polarization is not the first layer of a stratification. It is a Finsler Lie group.

We claim that the one-parameter subgroup $\sigma:\R\to G$, $\sigma(t):=\exp(tX_3)$ for all $t\in\R$, is a normal curve. Indeed, let $\lambda\in\g^*$ be the covector defined by the relations $\lambda(X_i)=0$ for $i\in\{1,2\}$,  $\lambda(X_3)=1$. It is straightforward to check that $\lambda\in \partial_{X_3}\E$. 
Since $X_3$ is in the center of $\g$, we have 
\begin{equation*}
    \lambda\Ad_{\exp(tX_3)}=\lambda\in \partial_{X_3}\E, \quad \forall t\in\R,
\end{equation*}
hence, the curve $\sigma$ solves the normal inclusion \eqref{eq:normal_equation}.
We claim that $\sigma$ is the unique normal curve with associated covector $\lambda$ starting at the origin.
Fix $x_1,x_2,x_3\in \R$. Notice that
\begin{align*}
\lambda\Ad_{\exp(x_1X_1+x_2X_2+x_3X_3)}(-x_2X_1+x_1X_2+X_3)&=x_1^2+x_2^2+1\\
&\geq \max\{|x_1|,|x_2|,1\}\\
&=\|(-x_2X_1+x_1X_2+X_3)\|.
\end{align*}
Consequently, we have that $
\|\lambda\Ad_{\exp(x_1X_1+x_2X_2+x_3X_3)}\|_*=1$ if and only if $x_1=x_2=0$. Hence, $\sigma$ is the unique normal curve with associated covector $\lambda$.

\begin{figure}[H]
\centering

\begin{tikzpicture}[scale=2]
\draw[->] (-0.5,0) -- (1.3,0) node[right] {$x$};
\draw[->] (0,-0.5) -- (0,1.3) node[above] {$y$};

\def\b{1} 

\draw[thick, blue, ->] (0,0) -- ({\b},0) node[midway, below]{$\beta$};

\draw[thick, blue, ->] ({\b},0) -- ({\b},{\b}) node[midway, right]{};

\draw[thick, blue, ->] ({\b},{\b}) -- (0,{\b}) node[midway, above]{};

\draw[thick, blue, ->] (0,{\b}) -- (0,0) node[midway, left]{$\beta$};

\end{tikzpicture}

\caption{The planar projection of the curve $\gamma$. It encloses area $\beta^2$ in the plane and it has length $4\beta$, like $\gamma$.}
\label{fig:gamma-xy-projection}

\end{figure}
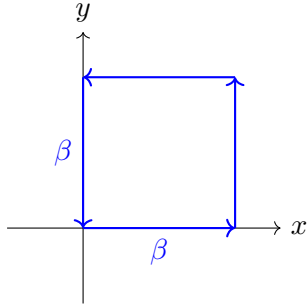

We next claim that the normal curve $\sigma$ is not locally optimal. 
Indeed, for each $\epsilon>0$ we will construct a curve that goes from $\exp(0)$ to $\exp(\epsilon X_3)$ and has length smaller than $\epsilon$. It will be made of the concatenation of four pieces of left translations of one-parameter subgroups.
Choose $\beta>0$ such that $4\beta+\beta^2=\epsilon$. 
Let $\gamma:[0,4\beta]\to \R^3$ be the curve with control:
\begin{equation*}
    u(t) := \begin{cases}
        X_1 + X_3 & \quad \text{for } t \in (0,\beta), \\
        X_2 + X_3 & \quad \text{for } t \in (\beta,2\beta), \\
        -X_1 + X_3 & \quad \text{for } t \in (2\beta,3\beta), \\
        -X_2 + X_3 & \quad \text{for } t \in (3\beta,4\beta). \\
    \end{cases}
\end{equation*}
We start by remarking that $\|u(t)\| = 1$ for a.e. $t \in [0,4\beta]$, hence the length of $\gamma$ is $4\beta < \epsilon$. A straightforward computation shows that
$\gamma(4\beta)
=\exp((4\beta+\beta^2) X_3)
=\exp(\epsilon X_3)$, see \cref{fig:gamma-xy-projection} .

\begin{remark}
 From the study of the one-parameter subgroup $\sigma$ above, we deduce that the Heisenberg group with the $L^\infty$ norm is an example of a Finsler Lie group that, as a metric space, is not locally uniquely geodesic. Indeed, notice that if in standard exponential coordinates, a curve $\kappa:=(\kappa_1,\kappa_2,\kappa_3)$ is a geodesic between $(0,0,0)$ and $(0,0,\epsilon)$, with $\epsilon>0$, then $(-\kappa_1,-\kappa_2,\kappa_3)$ is a geodesic between the same points, and it does not coincide with $\kappa$ because it is different than the one-parameter subgroup $t\mapsto \sigma(t)=(0,0,t)$.
\end{remark}

\bibliography{biblio}

\end{document}